\documentclass[10pt]{amsart}
\usepackage[dvips]{epsfig}
\usepackage{amsmath,amssymb,bbm}
\usepackage[all, knot]{xy}
\SelectTips{cm}{}
\xyoption{arc}
\xyoption{2cell}
\UseAllTwocells

\newtheorem{theorem}{Theorem}
\numberwithin{theorem}{subsection}
\newtheorem{proposition}[theorem]{Proposition}
\newtheorem{lemma}[theorem]{Lemma}
\newtheorem{varexample}[theorem]{Example}
\newenvironment{example}{\begin{varexample}\em}{\em\end{varexample}}
\newtheorem{definition}[theorem]{Definition}
\newtheorem{varremark}[theorem]{Remark}
\newenvironment{remark}{\begin{varremark}\em}{\em\end{varremark}}

\newcommand{\ra}{\mathop{\rightarrow}}
\newcommand{\ralim}{\mathop{\ra}\limits}
\newcommand{\la}{\mathop{\leftarrow}}
\newcommand{\lalim}{\mathop{\la}\limits}

\newcommand{\opname}[1]{\operatorname{#1}}

\newcommand{\catname}[1]{\boldsymbol{\opname{{#1}}}}
\newcommand{\inprod}[1]{\left\langle{#1}\right\rangle}
\newcommand{\Id}{\opname{Id}}

\newcommand{\C}{\catname{C}}

\newcommand{\X}{\catname{X}}
\newcommand{\Y}{\catname{Y}}

\newcommand{\Span}{\opname{Span}}

\newcommand{\V}{\catname{Vect}}
\newcommand{\iiV}{\catname{2Vect}}

\newcommand{\Set}{\catname{Set}}
\newcommand{\FinSet}{\catname{FinSet}}
\newcommand{\Gpd}{\catname{FinGpd}}
\newcommand{\Cat}{\catname{Cat}}
\newcommand{\HV}[1]{[#1,\V]}
\newcommand{\FV}{\Lambda}

\newcommand{\CG}{\mathbb{C}[G]}

\newcommand{\VG}{\catname{Vect}[G]}

\newcommand{\id}{\opname{id}}

\begin{document}

\title{2-Vector Spaces and Groupoids}
\author{Jeffrey C. Morton}
\address{Mathematics Department, University of Western Ontario}
\email{\tt{jeffrey.c.morton@gmail.com}}

\maketitle

\begin{abstract}
\addcontentsline{toc}{chapter}{Abstract}
This paper describes a relationship between essentially finite
groupoids and 2-vector spaces.  In particular, we show to construct
2-vector spaces of \textit{$\V$-valued presheaves} on such groupoids.
We define 2-linear maps corresponding to functors between groupoids in
both a covariant and contravariant way, which are ambidextrous
adjoints.  This is used to construct a representation---a weak
functor---from $\Span(\Gpd)$ (the bicategory of essentially finite groupoids and spans of
groupoids) into $\iiV$.  In this paper we prove this and give the
construction in detail.
\end{abstract}

\maketitle

\section{Introduction}\label{sec:intro}

In this paper, I will describe an extension of the
\textit{groupoidification} program of Baez and Dolan \cite{hdaVII}.
Groupoidification refers to the program of treating parts of linear
algebra as arising from spans of groupoids (categories whose morphisms
are all invertible) by a process of ``degroupoidification'', which
produces complex vector spaces associated to groupoids, and linear
maps associated to spans.  The extension described here shows a
connection of the setting of groupoids and spans with 2-vector spaces
and 2-linear maps, a categorical analog of linear algebra.  (We will
assume that all groupoids are essentially finite - that is, equivalent
to finite groupoids - although there is work in progress on how to
extend these results to infinite groupoids, and in particular Lie
groupoids.)

A simple example of the groupoidification program can be seen in terms
of spans of finite sets (i.e. finite trivial groupoids).  In that
program, groupoids give corresponding vector spaces and spans of
groupoids give corresponding linear maps.  In particular, the special
case of trivial groupoids (equivalently, sets) gives a useful
illustration.  Given a finite set $S$, there is a finite dimensional
vector space $L(S)$ consisting of all complex linear combinations of
elements of $S$.  Now, consider a span in $\FinSet$: that is, a
diagram of the form:
\begin{equation}
\xymatrix{
 & X \ar[dl]_{s} \ar[dr]^{t} & \\
Y & & Z
}
\end{equation}
To the span, there is a corresponding to a linear map $L(X) : L(Y)
\rightarrow L(Z)$, represented by a matrix $T$ whose $(i,j)$-component is
$|(s,t)^{-1}(y_i,z_j)$

\begin{figure}[h]
\begin{center}
\includegraphics{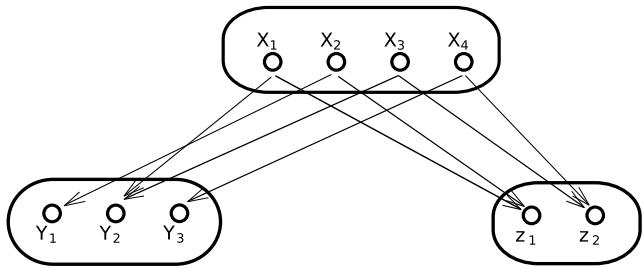}
\end{center}
\caption{\label{fig:set-span}A Span of Sets}
\end{figure}
So the set in Figure \ref{fig:set-span} gives rise to the linear
transformation:
\begin{align}
\nonumber y_1 \mapsto z_1 \\
y_2 \mapsto z_1 + z_2 \\
\nonumber y_3 \mapsto z_2
\end{align}

This makes sense for spans of finite sets.  Similarly, we will be
considering an analogous construction for spans of \textit{essentially
finite} groupoids.

There is a physical motivation here in quantum mechanics.  If $Y$ is
the (discrete) set of classical (pure) states of a system, then $L(Y) =
\mathbb{C}[Y]$, the space of linear combinations of states in $Y$, is
the Hilbert space of the corresponding quantum mechanical system.
(More generally, if $Y$ is a measure space, one takes $L^2(Y)$).
In the span, we think of $X$ as a set of ``processes'' $x$, each
with a designated ``source'' or starting state $s(x) \in Y$ and
``target'' or ending state $t(x) \in Z$.  Then the linear
transformation described by the matrix $T$ can be seen in the
following way, which we shall generalize later on:

Given a linear combination of elements of $Y$ (that is, a function $f
: Y \rightarrow \mathbb{C}$), transport $f$ to $X$ by ``pulling
back'' along $s$.  That is, $s^{\ast}f(x) = f(s(x))$.  Then ``push
forward'' to $Z$ by taking the sum over all elements of $X$ mapping
down to a chosen one in $Z$:
\begin{equation}
t_{\ast} s^{\ast} f (z) = \sum_{t(x) = z} s^{\ast}f (x)
\end{equation}
This precisely gives matrix multiplication by the matrix described
above, and can clearly also be seen as a ``sum over histories'': the
value of $t_{\ast} s^{\ast} f(z)$ is a sum, over all histories $x$
ending at $z$, of the value of $f$ at the source $s(x)$.  This
illustrates a contrast between classical and quantum
\textit{processes}.  Classically, states succeed each other by exactly
one process.  In the quantum picture, \textit{every} possible process
contributes to evolution of a state.  In particular, there is an
interpretation of quantum processes in terms of ``matrix mechanics'',
which takes a sum (in the form of matrix multiplication) over all
histories joining fixed start and end states.  This is exactly what is
shown in our example.

It is not too difficult to check that the linearization of spans of
sets gets along with composition, so that the composite of spans (by
pullback, giving a set of \textit{composite paths}) agrees with
composition of linear maps.  That is, that the process is
\textit{functorial}.  This fact makes it possible to think of
categorifying this process, in order to explicitly include symmetries
of both states and histories as fundamental concepts.  A categorified
version of this process should be a 2-functor.

One way to generalize spans of sets, which is seen in \cite{hdaVII},
uses groupoids (categories whose morphisms are all invertible) instead
of sets.  One reason to consider this is that it often happens that
the configuration space can naturally be thought of not as a set but
as a groupoid.  This happens particularly when there are symmetry
operations acting on the set of configurations, and we explicitly
represent such symmetries as morphisms of the groupoid.  The existence
of a group action on the set would be one example.  In such a
categorified picture, $X$ has objects which represent states of a
system, and morphisms denoting \textit{symmetries} of states.  Then
$L$ gives vector spaces which are linear combinations of
\textit{isomorphism classes of objects} of the groupoids.  The
components of the linear maps uses groupoid cardinality instead of set
cardinality:
\begin{align}
\label{eq:baezgpdify} L(X)_{[y_i],[z_k]} & = \sum_{x \in \underline{(s,t)^{-1}(y_i,z_k)}}\frac{\#(\opname{Aut}(y_i))}{\#(\opname{Aut}(x))} \\ 
\nonumber & = |\widehat{(y_i,z_k)}| \cdot \# (\opname{Aut}(y_i))
\end{align}
where $\widehat{(y_i,z_k)}$ is the \textit{essential preimage} of
$y_i$ and $z_k$, and its cardinality is the groupoid cardinality
described by Baez and Dolan \cite{finfeyn} (the other cardinality is
the order of the group).  This uses a weighting of contributions from
intermediate elements depending on the size of their symmetry group.
The groupoid cardinality of a finite groupoid $X$ is:
\begin{equation}\label{eq:gpdcard}
|X| = \sum_{[x] \in \underline{X}} \frac{1}{\# \opname{Aut}(x)}
\end{equation}
where the cardinality in the sum denotes the order of the group.

Here, however, we want to do something a little different: this
process is still a functor, and we wanted a \textit{2-functor}.  Since
we want to think of $X$ as a category, rather than look at functions
from the objects of $X$ into $\mathbb{C}$, we should look at functors
from $X$ into some category which plays the role of $\mathbb{C}$.  In
particular, this category will be $\V$, whose objects are vector
spaces over $\mathbb{C}$, and whose morphisms are linear maps.  When
categorifying, therefore, we will want to find an analogous 2-functor,
which requires specifying more data.

Then there will be a ``free 2-vector space'' $\FV(X)$ of all functors
from $X$ into $\V$.  We think of the objects as ``2-linear
combinations'' of classical states, each with an \textit{internal
state space} which carries a representation of the symmetry group of
that state.  For most physically realistic systems, $X$ would be an
infinite set with a measure, and in fact a symplectic manifold.  In
general, to deal with $L^2$ spaces involves some issues in analysis,
such as the measure on $X$.  Then instead of $L(X)$ we consider
$L^2(X)$.  A similar caveat should apply in the categorified setting.
Restricting to the situation of a finite groupoid helps to more
clearly illustrate some of the purely category-theoretic aspects of
the ``free 2-vector space'' construction.  We do expect that for
well-behaved smooth groupoids, for example, similar results to those
considered in this paper will hold, involving infinite dimensional
2-vector spaces one could denote $2L^2(X)$.  But this will be
addressed in a companion paper.

Finally, we remark that this construction is used in the construction
of an Extended Topological Quantum Field Theory (ETQFT) in the
author's Ph.D. thesis \cite{mortonthesis}, where the groupoids in
question are topological invariants of manifolds.  By analogy, it
could be used to give ``extended quantum theories'' in other settings
where spans of groupoids appear.

Another view of a related process involves the 2-functor into additive
categories (which, in the $\mathbb{C}$-linear case, are the KV
2-vector spaces) from a 2-category $\catname{Bim}$, whose objects are
rings, morphisms in $\hom_{\catname{Bim}}(R,S)$ are $(R,S)$-bimodules,
and whose 2-morphisms are bimodule homomorphisms.  This is a dual
picture to that of spans.  Indeed, the type of ``pull-push''
construction given here is ubiquitous (as its appearance in linear
algebra suggests), due to the universal properties of categories of
spans (see \cite{unispan}).  Similar notions also appear in the theory
of Mackey functors (see \cite{elango}), and in the study of
``correspondences'' in noncommutative geometry, algebraic geometry,
and elsewhere.

In this paper, we will begin by describing the source and target
categories, $\Span(\Gpd)$ in section \ref{sec:spangpd}, and $\iiV$ in
section \ref{sec:KV}.  In particular, to categorify the functor $L$,
we need a 2-category to correspond to $\V$, and this will be the
2-category of all Kapranov-Voevodsky 2-vector spaces.  A KV 2-vector
space is an abelian category with some extra structure, just as a
vector space is a special type of abelian group.  In section
\ref{sec:KV} we give some background and collect some fundamental
results about them which are widely known, but whose proofs are seldom
given.  For example, we show that 2-vector spaces, understood as a
semisimple $\mathbb{C}$-linear additive category, are all equivalent
to $\V^k$ for some nonnegative integer $k$.

In section \ref{sec:kv2vsgrpd}, we give the object level of our
construction for a 2-functor $\FV$, which, to (essentially finite)
groupoids assigns KV 2-vector spaces.  Analogously with sets, we
obtain 2-linear maps for spans of groupoids.  In fact, just as with
sets, this is a consequence of an even simpler correspondence.
Namely, there is the ``pullback'' and ``push-forward'' of a function
mentioned in the description of the linear map from a span of sets as
a \textit{sum over histories} (the sum occurs in the ``push-forward''
operation, and corresponds to the sum in matrix multiplication).  The
groupoid situation is more complicated than that for sets, however,
because of the existence of automorphisms of the objects, and the
condition that maps between groupoids are functors.  This means, in
particular, that for each object $x$ in a groupoid, the functor
determines a homomorphism from the automorphism group of $x$ to that
of its image.  The push-forward operation can be interpreted as a Kan
extension and has both an object and a morphism level.

Section \ref{sec:spans} describes how these results define 2-linear
maps associated to spans of groupoids.  It begins with a brief
discussion of the bicategory whose objects are groupoids, whose
morphisms are spans of groupoids, and whose 2-morphisms are span maps.
This is followed by an explicit construction of the morphism level of
the 2-functor $\FV : \Span(\Gpd) \rightarrow \iiV$ and shows that it
preserves composition of spans in the weak sense: that is, up to a
specified isomorphism. (Technical details of this proof are reserved
for appendix \ref{ax:composproof}).  Finally, using Frobenius
reciprocity, it describes a simple explicit matrix representation for
the 2-functor constructed.

Section \ref{sec:spanspan} then continues by describing how this
representation works at the level of 2-morphisms.  This is analogous
to the 1-morphism level, in that it consists of a ``pullback and
pushforward'' process.  This is most easily described in terms of the
linear maps between corresponding vector spaces which appear in the
matrix representation of the 2-linear maps associated to a pair of
spans from $A$ to $B$.  We give this construction and show it
preserves both vertical and horizontal composition of 2-morphisms
in the appropriate ways.

The results shown in sections \ref{sec:spans} and \ref{sec:spanspan}
do much of the work involved in showing that our representation is
really a 2-functor.  The remainder of this proof is given in section
\ref{sec:maintheorem}.

Now we begin to describe the 2-linearization process by collecting
some key facts about the bicategory of 2-vector spaces, including a
canonical construction of one for each (essentially finite) groupoid.

\section{The Bicategory $\Span(\Gpd)$}\label{sec:spangpd}

The main purpose of this paper is to describe a weak 2-functor
\begin{equation}
\FV : \Span(\Gpd) \ra \iiV
\end{equation}
In this section, we will describe the source bicategory, $\Span(\Gpd)$.

First, the objects of $\Span(\Gpd)$.

\begin{definition}\label{def:essfingpd} An \textbf{essentially finite}
groupoid is one which is equivalent to a finite groupoid.  A
\textbf{finitely generated} groupoid is one with a finite set of
objects, and all of whose morphisms are generated under composition by
a finite set of morphisms.  An \textbf{essentially finitely generated}
groupoid is one which is equivalent to a finitely generated one.
\end{definition}

Note that, in particular, essentially finitely generated groupoids must be essentially finite, since every object has an identity morphism.  We will use the term ``essentially finite'' to mean both of these conditions.  Next we describe the morphisms of $\Span(\Gpd)$.

In any category $\C$, a \textit{span} is a diagram of the form:
\begin{equation}\label{eq:spandiag}
\xymatrix{
  & X \ar[dl]^{s} \ar[dr]_{t} & \\
A_1  & & A_2  \\
}
\end{equation}

\begin{figure}[h]
\begin{center}
\includegraphics{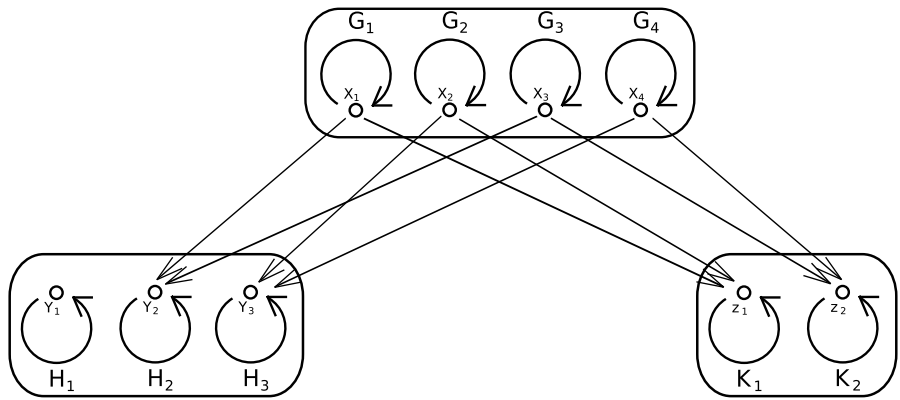}
\end{center}
\caption{\label{fig:grpd-span}A Span of Groupoids}
\end{figure}

In particular, we want to reproduce the ``linearization'' associated
to spans of sets which we discussed in the introduction.  The idea is
that given a span of groupoids, as in Figure \ref{fig:grpd-span}
(which suppresses the homomorphisms labelling the strands in the span,
but should be compared with Figure \ref{fig:set-span}), there will be
a ``transfer'' 2-linear map from the KV 2-vector space associated to
the source of the span, to that associated to the target.

If $\C$ has pullbacks, we can define composition of spans using them:
\begin{equation}\label{xy:spancomposite}
  \xymatrix{
     & & X' \circ X \ar[dl]_{S} \ar[dr]^{T} & & \\
     & X \ar[dl]_{s} \ar[dr]^{t} & & X' \ar[dl]_{s'} \ar[dr]^{t'} & \\ 
   A_1 & & A_2 & & A_3 \\ 
  } 
\end{equation} 
where we define $X' \circ X$ to be the object, unique up to
isomorphism, which makes the central square a pullback square.  That
is, it is a terminal occupant of this niche.  If $\C$ is, in addition,
a concrete category, the pullback is a subobject of the product $X
\times X'$.
\begin{equation}
X' \circ X = X \times_{A_2} X' = \bigcup_{a \in A_2} t^{-1}(a) \times (s')^{-1}(a)
\end{equation}
the fibred product of $X'$ and $X$ over $A_2$.  (Indeed, if $\C$ is Cartesian, any span can be factored
through a product.)


Now, for any category $\C$ with pullbacks, there is a category
$\Span(\C)$ whose objects are the objects of $\C$, and whose morphisms
are isomorphism classes spans in $\C$ composed by pullback.  Here we
are taking spans up to isomorphisms $\alpha : X_1 \rightarrow X_2$
which are commuting diagrams of the form:
\begin{equation}\label{xy:spanmap}
\xymatrix{
  & X_1 \ar[dl]_{s_1} \ar[d]^{\alpha} \ar[dr]^{t_1} & \\
A_1 & X_2 \ar[l]^{s_2} \ar[r]_{t_2} & A_2 \\
}
\end{equation}

However, we will do something slightly different.  We will be
interested in spans of groupoids.  Since groupoids naturally form a
2-category, we should weaken the notion of composition, and give an
appropriate notion of 2-morphism, for a bicategory $\Span(\Gpd)$.

Given any 2-category $\C$ with weak pullbacks, one can again form a
bicategory $\Span(\C)$ with the same objects as $\C$ and with spans
for morphisms.  T composite of spans $A_1 \lalim^{s} X \ralim^{t} A_2$
and $A_2 \lalim^{s'} X' \ralim^{t'} A_3$ is an object $X' \circ X$
together with a 2-morphism $\alpha$ making this diagram commute, and
terminal for such choices:
\begin{equation}\label{xy:wkpb}
 \xymatrix{
   & & {X'\circ X} \ar_{S}[ld]\ar^{T}[rd]\ar@/_2pc/_{s \circ S}[ddll]\ar@/^2pc/^{t' \circ T}[ddrr] & & \\
   & {X}\ar^{s}[ld]\ar_{t}[rd] \ar@{=>}[rr]^{\alpha}_{\sim} & & {X'}\ar^{s'}[ld]\ar_{t'}[rd] & \\
   {A_1} & & {A_2} & & {A_3} \\
}
\end{equation}

A 2-morphism in $\Span(\C)$ will be an isomorphism class of
\textit{spans of span maps}.  That is, consider a span of
\textit{2-morphisms} in the usual $\Span(\C)$.  This is a diagram of
the form:
\begin{equation}\label{xy:spanmapspan}
\xymatrix{
  & X_1 \ar[dl]^{s_1}  \ar[dr]_{t_1} & \\
A_1 & Y \ar[l] \ar[r] \ar[d]^{t} \ar[u]_{s}& A_2 \\
  & X_2 \ar[ul]_{s_2}  \ar[ur]^{t_2} & \\
}
\end{equation}

In principle we need only require this diagram commute weakly: that
is, there are isomorphisms $\zeta_s : s_1 \circ s \ra s_2 \circ t$ and
$\zeta_t : t_1 \circ s \ra t_2 \circ t$.  For the most part, since the
construction we mean to give is invariant under equivalence of
groupoids, and taking the $A_i$ to be skeletal makes this strict, we
will assume strict commutativity, though we shall indicate where the
argument must be changed to accommodate the weak case.

We are considering such diagrams only \textit{up to isomorphism}: that
is, the inner span $X_1 \la Y \ra X_2$ in the 2-morphism
\ref{xy:spanmapspan} is only considered up to an isomorphism of spans
in the sense of \ref{xy:spanmap}.

The reason for considering these is as follows.  First, taking a
category $\C$ and passing to $\Span(\C)$ amounts to formally adjoining
duals for morphisms in $\C$.  The dual of any span is the span which
has the same maps, considered in the reverse orientation, exchanging
the role of source and target object.  When we apply $\FV$, these
duals will in fact become adjoints, as we shall see.

Now, we are interested in the case $\C = \Gpd$, so the diagram for the
composite of spans between groupoids contains a weak pullback square:
composition is only preserved up to isomorphism.  In particular, the
objects are now groupoids, which are themselves categories with
objects and morphisms.  Since it makes sense to speak of two objects
of a groupoid being isomorphic, the weakest meaningful condition is
that objects of groupoids $X$ and $X'$ should need only project to
isomorphic objects on $A_2$.  But there are potentially different
isomorphisms between those objects.  So the weak pullback is a larger
groupoid than a strict pullback, since its objects come with a
specified isomorphism between the two restrictions.

That this is a weak pullback square of functors between groupoids
means that this diagram commutes up to the natural isomorphism $\alpha
: t \circ S \longrightarrow s' \circ T$.  (The fact that $\alpha$ is
iso is what makes this a \textit{weak} pullback rather than a
\textit{lax} pullback, where $\alpha$ is only a natural
transformation.)  This is an example of a \textit{comma category} (the
concept, though not the name, was introduced by Lawvere in his
doctoral thesis \cite{lawverethesis}).  We recall some background
about this construction in Appendix \ref{sec:commacat}.

Now, the process of finding higher morphisms by taking spans of span
maps could obviously be continued: each new level of span naturally
gives maps of spans as morphisms.  We could repeat the process of
adjoining duals by passing to spans of such span maps, and so on
recursively as far as we wish.  For our purposes here, however, we
will stop at 2-morphisms for two reasons.  First, we want to describe
a representation into $\iiV$, which is a 2-category.  This in turn is
since the objects of $\Gpd$ are themselves categories, and our
2-functor $\FV$ will represent them as categories - so a 2-category is
the natural setting for them.

Collecting the definition together, we then have the following.

\begin{definition}\label{def:spangpd}
The bicategory $\Span(\Gpd)$ has:
\begin{itemize}
\item \textbf{Objects}: Essentially finite groupoids
\item \textbf{Morphisms}: Spans of groupoids, composed by weak pullback
\item \textbf{2-Morphisms}: Isomorphism classes of spans of span maps,
  composed by weak pullback both horizontally and vertically
\end{itemize}
\end{definition}

Having defined $\Span(\Gpd)$, the source bicategory of the 2-functor $\FV$ we aim to describe here, we next describe its target, $\iiV$. 

\section{Kapranov-Voevodsky 2-Vector Spaces}\label{sec:KV}

There are two major philosophies regarding how to categorify the
concept ``vector space''.  A Baez-Crans (BC) 2-vector space is a
category object in $\V$---that is, a category having a vector space of
objects and of morphisms, where source, target, composition, etc. are
linear maps.  This is a useful concept for some purposes---it was
developed to give a categorification of Lie algebras.  The reader may
refer to the paper of Baez and Crans \cite{BC} for more details.
However, a BC 2-vector space turns out to be equivalent to a 2-term
chain complex and for many purposes this is too strict.  This is not
the concept of 2-vector space which concerns us here.

The other, earlier, approach is to define a 2-vector space as a
category having operations such as a monoidal structure analogous to
the addition on a vector space.  In particular, we will restrict our
attention to \textit{complex} 2-vector spaces.

This ambiguity about the correct notion of ``2-vector space'' is
typical of the problem of categorification.  Since the categorified
setting has more layers of structure, there is a choice of level to
which the structure in the concept of a vector space should be lifted.
Thus in the BC 2-vector spaces, we have literal vector addition and
scalar multiplication within the objects and morphisms.  In KV
2-vector spaces and their cousins, we only have this for morphisms,
and for objects there is a categorified analog of addition, in the
sense that they are additive categories.  The key difference between
the two notions of 2-vector space lies in which category plays the
role of the ``base field'': in the BC definition, this is the ring
category $\mathbb{C}_{[0]}$ whose objects are complex numbers, whereas
for the KV definition it is $\V$, whose objects are complex vector
spaces.  This is discussed by Josep Elgueta \cite{elgueta}.

Indeed, Elgueta \cite{elgueta} shows several different types of
``generalized'' 2-vector spaces, and relationships among them.  In
particular, while KV 2-vector spaces can be thought of as having a
\textit{set} of basis elements, a generalized 2-vector space may have
a general \textit{category} of basis elements.  The free generalized
2-vector space on a category is denoted $\catname{Vect[\mathcal{C}]}$.
Then KV 2-vector spaces arise when $\mathcal{C}$ is a discrete
category with only identity morphisms.  This is essentially a set $S$
of objects.  Thus it should not be surprising that KV 2-vector spaces
have a structure analogous to free vector spaces generated by some
finite set - which are isomorphic to $\mathbb{C}^k$.

\subsection{Definition}

The standard example of this approach is the Kapranov-Voevodsky (KV)
definition of a 2-vector space \cite{KV}, which is the form we shall
use (at least when the situation is finite-dimensional).  To motivate
the KV definition, consider the idea that, in categorifying, one
should replace the base field $\mathbb{C}$ with a monoidal category.
Specifically, it turns out, with $\V$, the category of finite
dimensional complex vector spaces.  This leads to the following
replacements for concepts in elementary linear algebra:
\begin{itemize}
\item {Vectors = $k$-tuples of scalars} $\mapsto$ {2-vectors = $k$-tuples of vector spaces}
\item {Addition} $\mapsto$ {Direct Sum}
\item {Multiplication} $\mapsto$ {Tensor Product}
\end{itemize} So just as $\mathbb{C}^k$ is the standard example of a
complex vector space, $\V^k$ will be the standard example of a
2-vector space.  But we should define these precisely.

To begin with, a KV 2-vector space is a $\mathbb{C}$-linear additive
category with some properties, so we begin by explaining this.  The
property of \textit{additivity} for categories, is here seen as the
analog of the group structure of a vector space, though additivity in
a category is somewhat different.  The motivating example for us is
the \textit{direct sum} operation in $\V$.  Such an operation plays
the role in a 2-vector space which vector addition plays in a vector
space.

\begin{definition}
If a category $\C$ is enriched in abelian groups, a \textbf{biproduct} is an operation giving, for
any objects $x$ and $y$ in $\C$ an object $x \oplus y$
equipped with morphisms $\iota_x,\iota_y$ from $x$ and $y$
respectively into $x \oplus y$; and morphisms $\pi_x, \pi_y$ from $x
\oplus y$ into $x$ and $y$ respectively, which satisfy the
biproduct relations:
\begin{equation}\label{eq:biproduct1}
\pi_x \circ \iota_x = \id_x \text{  and  } \pi_y \circ \iota_y = \id_y
\end{equation}
and
\begin{equation}\label{eq:biproduct2}
\iota_x \circ \pi_x + \iota_y \circ \pi_y = \id_{x \oplus y}
\end{equation}
\end{definition}
Whenever biproducts exist, they are always both products and
coproducts.

\begin{definition} A \textbf{$\mathbb{C}$-linear additive category}
is a category $\catname{V}$ enriched in $\V$(i.e. $\forall x, y \in \catname{V}$, $\hom(x,y)$ is a vector space
over $\mathbb{C}$), such that composition is a bilinear map, and such that $\catname{V}$ has a zero object (i.e. $0$ which is both initial and terminal).  A
$\mathbb{C}$-linear functor between $\mathbb{C}$-linear categories
is one where morphism maps are $\mathbb{C}$-linear.  A \textbf{simple object} in $\catname{V}$ is $x \in \catname{V}$
such that $\hom(x,x) \cong \mathbb{C}$.
\end{definition}

As important fact about KV 2-vector spaces is that they have
(finite) bases: they are generated by finitely many simple objects.

\begin{definition} A \textbf{Kapranov--Voevodsky 2-vector space} is a
$\mathbb{C}$-linear additive category which is semisimple (every object can be
written as a finite biproduct of simple objects).  A \textbf{2-linear map} between
2-vector spaces is a $\mathbb{C}$-linear functor.
\end{definition}

\begin{remark}It is a consequence of $\mathbb{C}$-linearity that a 2-linear map also preserves biproducts, since the images of the $\pi$ and $\iota$ maps still satisfy the definition of a biproduct (and the universal properties for product and coproduct follow automatically).  The above definition of a
2-linear map is sometimes given in the equivalent form requiring that
the functor preserve exact sequences.  Indeed, since every object is a
finite biproduct of simple objects, a 2-vector space is an
\textit{abelian} category.  (See e.g. Freyd \cite{freyd}.)
\end{remark}



\begin{example}The standard example \cite{KV} of a KV 2-vector space
highlights the analogy with the familiar vector space $\mathbb{C}^k$.
The 2-vector space $\V^k$ is a category whose objects are $k$-tuples
of vector spaces, maps are $k$-tuples of linear maps.  The
\textit{additive} structure of the 2-vector space $\V^k$ comes from
applying the direct sum in $\V$ component-wise.

Note that there is an equivalent of \textit{scalar multiplication},
using the tensor product:
\begin{equation}
V \otimes 
\begin{pmatrix}
  V_1 \\ 
  \vdots \\
  V_k \\
\end{pmatrix}
=
\begin{pmatrix}
  V \otimes V_1 \\
  \vdots \\
  V \otimes V_k
\end{pmatrix}
\end{equation}
and
\begin{equation}
\begin{pmatrix}
  V_1 \\ 
  \vdots \\
  V_k
\end{pmatrix}
 \oplus 
\begin{pmatrix}
  W_1 \\
  \vdots \\
  W_k
\end{pmatrix}
= 
\begin{pmatrix}
  V_1 \oplus W_1 \\
  \vdots \\
  V_k \oplus W_k
\end{pmatrix}
\end{equation}

As the correspondence with linear algebra would suggest, 2-linear maps
$T: \catname{Vect^k} \ra\catname{Vect^l}$ amount to $k \times l$
matrices of vector spaces, acting by matrix multiplication using the
direct sum and tensor product instead of operations in $\mathbb{C}$:
\begin{equation}\label{eq:kv2linmatrix}
\begin{pmatrix}
T_{1,1} & \dots & T_{1,k} \\
\vdots & & \vdots \\
T_{l,1} & \dots & T_{l,k} \\
\end{pmatrix}
\begin{pmatrix}
  V_1 \\ 
  \vdots \\
  V_k
\end{pmatrix}
=
\begin{pmatrix}
  \bigoplus_{i=1}^k T_{1,i} \otimes V_i \\ 
  \vdots \\
  \bigoplus_{i=1}^k T_{l,i} \otimes V_i \\ 
\end{pmatrix}
\end{equation}

The natural transformations between these are matrices of linear
transformations:
\begin{equation}\label{eq:kvnattransmatrix}
\alpha = \begin{pmatrix}
\alpha_{1,1} & \dots & \alpha_{1,k} \\
\vdots & & \vdots \\
\alpha_{l,1} & \dots & \alpha_{l,k} \\
\end{pmatrix}
:
\begin{pmatrix}
T_{1,1} & \dots & T_{1,k} \\
\vdots & & \vdots \\
T_{l,1} & \dots & T_{l,k} \\
\end{pmatrix}
\longrightarrow
\begin{pmatrix}
T'_{1,1} & \dots & T'_{1,k} \\
\vdots & & \vdots \\
T'_{l,1} & \dots & T'_{l,k} \\
\end{pmatrix}
\end{equation}
where each $\alpha_{i,j} : T_{i,j} \ra T'_{i,j}$ is a linear
map in the usual sense.

These natural transformations give 2-morphisms between 2-linear maps, so
that $\V^k$ is a bicategory with these as 2-cells:
  \begin{equation}\label{xy:2vs2cell}
    \xymatrix{
      \catname{Vect^k} \ar@/^1pc/[r]^{F}="0"
      \ar@/_1pc/[r]_{G}="1" & \catname{Vect^l} \\ \ar@{=>}"0"
      ;"1"^{\alpha}
    }
  \end{equation}

In our example above, the finite set of simple objects of which every
object is a sum is the set of 2-vectors of the form
\begin{equation}
\begin{pmatrix}
0 \\
\vdots \\
\mathbb{C} \\
\vdots \\
0
\end{pmatrix}
\end{equation}
which have the zero vector space in all components except one (which
can be arbitrary).  We can call these \textit{standard basis
2-vectors}.  Clearly every object of $\V^k$ is a finite biproduct of
these objects, and each is simple (its vector space of endomorphisms
is 1-dimensional).
\end{example}

\subsection{Classification Theorems}

The most immediately useful fact about KV 2-vector spaces is the
following well known characterization:

\begin{theorem}\label{thm:kvvn} Every KV 2-vector space is equivalent as a category
to $\V^k$ for some $k \in \mathbb{N}$.
\end{theorem}
\begin{proof}
Suppose $\catname{K}$ is a KV 2-vector space with a basis of simple
objects $X_1 \dots X_k$.  Then we construct an equivalence $E:
\catname{K} \ra \V^k$ as follows:

$E$ should be an additive functor with $E(X_i) = V_i$, where $V_i$ is
the $k$-tuple of vector spaces having the zero vector space in every
position except the $i^{th}$, which has a copy of $\mathbb{C}$.  But
any object $X$, is a sum $\bigoplus_{i} X_i^{n_i}$, so by linearity
(i.e. the fact that $E$ preserves biproducts) $X$ will be sent to the
sum of the same number of copies of the $V_i$, which is just a
$k$-tuple of vector spaces whose $i^{th}$ component is
$\mathbb{C}^{n_i}$.  So every object in $K$ is sent to an $k$-tuple
of vector spaces.  By $\mathbb{C}$-linearity, and the fact that
hom-vector spaces of simple objects are one-dimensional, this
determines the images of all morphisms.

But then the weak inverse of $E$ is easy to construct, since sending
$V_i$ to $X_i$ gives an inverse at the level of objects, by the same
linearity argument as above. At the level of morphisms, the same
argument holds again.
\end{proof}

This is a higher analog of the fact that every finite dimensional
complex vector space is isomorphic to $\mathbb{C}^k$ for some $k \in
\mathbb{N}$.  So, indeed, the characterization of 2-vector spaces in
our example above is generic: every KV 2-vector space is equivalent to
one of the form given.  Moreover, our picture of 2-linear maps is also
generic, as shown by this argument, analogous to the linear algebra
argument for representation of linear maps by matrices:

\begin{lemma}\label{lemma:kv2linmatrix} Any 2-linear map $T: \V^n \ra
\V^m$ is naturally isomorphic to a map of the form (\ref{eq:kv2linmatrix}).
\end{lemma}
\begin{proof}
Any 2-linear map $T$ is a $\mathbb{C}$-linear additive functor
between 2-vector spaces.  Since any object in a 2-vector space can be
represented as a biproduct of simple objects---and morphisms
likewise---such a functor is completely determined by its effect on
the basis of simple objects and morphisms between them.

But then note that since the automorphism group of a simple object is
by definition just all (complex) multiples of the identity morphism,
there is no choice about where to send any such morphism.  So a
functor is completely determined by the images of the basis objects,
namely the 2-vectors $V_i = (0,\dots,\mathbb{C},\dots,0) \in \V^n$,
where $V_i$ has only the $i^{th}$ entry non-zero.

On the other hand, for any $i$, $T(V_i)$ is a direct sum of some
simple objects in $\V^m$, which is just some 2-vector, namely a
$k$-tuple of vector spaces.  Then the fact that the functor is
additive means that it has exactly the form given.
\end{proof}

And finally, the analogous fact holds for natural transformations
between 2-linear maps:

\begin{lemma}\label{lemma:kvnattransmatrix}Any natural transformation
$\alpha : T \ra T'$ from a 2-linear map $T : \V^n \ra \V^m$ to a
2-linear map $T' : \V^n \ra \V^m$, both in the form
(\ref{eq:kv2linmatrix}) is of the form (\ref{eq:kvnattransmatrix}).
\end{lemma}
\begin{proof}
By Lemma \ref{lemma:kv2linmatrix}, the 2-linear maps $T$ and $T'$ can
be represented as matrices of vector spaces, which act on an object
in $\V^n$ as in (\ref{eq:kv2linmatrix}).  A natural transformation
$\alpha$ between these should assign, to every object $X \in \V^n$, a
morphism $\alpha_X : T(X) \ra T'(X) $ in $\V^m$, such that the usual
naturality square commutes for every morphism $f : X \ra Y$ in $\V^n$.

Suppose $X$ is the $n$-tuple $(X_1, \dots , X_n)$,
where the $X_i$ are finite dimensional vector spaces.  Then
\begin{equation}\label{eq:TX}
T(X) = (\oplus_{k=1}^{n} V_{1,k} \otimes X_k, \dots,\oplus_{k=1}^{n} V_{m,k} \otimes X_k )
\end{equation}
where the $V_{i,j}$ are the components of $T$, and similarly
\begin{equation}
T'(X) = (\oplus_{k=1}^{n} V'_{1,k} \otimes X_k, \dots,\oplus_{k=1}^{n} V'_{m,k} \otimes X_k )\end{equation}
where the $V'_{i,j}$ are the components of $T'$.

Then a morphism $\alpha_X : T(X) \ra T'(X)$ consists of an $m$-tuple
of linear maps:
\begin{equation}
\alpha_j : \oplus_{k=1}^{n} V_{j,k} \otimes X_k \ra \oplus_{k=1}^{n} V'_{j,k} \otimes X_k
\end{equation}
but by the universal property of the biproduct, this is the same as
having an $(n \times m)$-indexed set of maps
\begin{equation}
\alpha_{jk} : V_{j,k} \otimes X_k \ra \oplus_{r=1}^{n} V'_{j,r} \otimes X_r
\end{equation}
and by the dual universal property, this is the same as having $(n
\times n \times m)$-indexed maps
\begin{equation}
\alpha_{jkr} : V_{j,k} \otimes X_k \ra V'_{j,r} \otimes X_r
\end{equation}
However, we must have the naturality condition for every morphism
$f:X \ra X'$:
\begin{equation}
\xymatrix{
  T(X) \ar[d]_{\alpha_{X}} \ar[r]^{T(f)} & T(X') \ar[d]^{\alpha_{X'}} \\
  T'(X) \ar[r]_{T'(f)} & T'(X')\\ 
}
\end{equation}
Note that each of the arrows in this diagram is a morphism in $\V^m$,
which are linear maps in each component---so in fact we have a
separate naturality square for each component.

Also, since $T$ and $T'$ act on $X$ and $X'$ by tensoring
with fixed vector spaces as in (\ref{eq:TX}), one has $T(f)_i =
\oplus_j 1_{V_{ij}} \otimes f_j$, having no effect on the $V_{ij}$.
We want to show that the components of $\alpha$ affect \textit{only}
the $V_{ij}$.

Additivity of all the functors involved implies that the assignment
$\alpha$ of maps to objects in $\V^n$ is additive.  So consider the
case when $X$ is one of the standard basis 2-vectors, having
$\mathbb{C}$ in one position (say, the $k^{th}$), and the zero vector
space in every other position.  Then, restricting to the naturality
square in the $k^{th}$ position, the above condition amounts to having
$m$ maps (indexed by $j$):
\begin{equation}
\alpha_{j,k} : V_{j,k} \ra V'_{j,k}
\end{equation}
So by linearity, a natural transformation is determined by an $n\times
m$ matrix of maps as in (\ref{eq:kvnattransmatrix}).
\end{proof}

The fact that 2-linear maps between 2-vector spaces are functors
between categories recalls the analogy between linear algebra and
category theory in the concept of an \textit{adjoint}.  If $V$ and $W$
are inner product spaces, the adjoint of a linear map $F : V \ra W$ is
a map $F^{\dagger}$ for which $\inprod{Fx,y} =
\inprod{x,F^{\dagger}y}$ for all $x \in V_1$ and $y \in V_2$.  A
(right) adjoint of a functor $F : \catname{C} \ra \catname{D}$ is a
functor $G : \catname{D} \ra \catname{C}$ for which $\hom_D(Fx,y)
\cong \hom_C(x,Gy)$ (and then $F$ is a left adjoint of $G$).

In the situation of a KV 2-vector space, the categorified analog of
the adjoint of a linear map is indeed an adjoint functor.  (Note that
since a KV 2-vector space has a specified basis of simple objects, it
makes sense to compare it to an inner product space.)  Moreover, the
adjoint of a functor has a matrix representation which is much like
the matrix representation of the adjoint of a linear map.  We
summarize this in the following (a variant of proposition 25 in \cite{hdaII}, shown there for 2-Hilbert spaces):

\begin{theorem}\label{thm:biadjoint} Given any 2-linear map $F : V \ra
W$, there is a 2-linear map $F^{\dagger} : W \ra V$ which is both a left and
right adjoint to $F$.
\end{theorem}
\begin{proof}
By Theorem \ref{thm:kvvn}, we have $V \simeq \V^n$ and $W \simeq \V^m$
for some $n$ and $m$.  By composition with these equivalences, we can
restrict to this case.  But then we have by Lemma
\ref{lemma:kv2linmatrix} that $F$ is naturally isomorphic to some
2-linear map given by matrix multiplication by some matrix of vector
spaces $[F_{i,j}]$:
\begin{equation}
\begin{pmatrix}
F_{1,1} & \dots & F_{1,n} \\
\vdots & & \vdots \\
F_{m,1} & \dots & F_{m,n} \\
\end{pmatrix}
\end{equation}

We claim that a (two-sided) adjoint functor $F^{\dagger}$ is given by the ``dual
transpose matrix'' of vector spaces $[F_{i,j}]^{\dagger}$:
\begin{equation}
\begin{pmatrix}
F^{\dagger}_{1,1} & \dots & F^{\dagger}_{1,m}\\
\vdots & & \vdots \\
F^{\dagger}_{n,1} & \dots & F^{\dagger}_{n,m} \\
\end{pmatrix}
\end{equation}
where $F^{\dagger}_{i,j}$ is the vector space dual $(F_{j,i})^{\ast}$
(note the transposition of the matrix).

We note that this prescription is symmetric, since $[T]^{\dagger
\dagger} = [T]$, so if $F^{\dagger}$ is always a left adjoint of $F$, then $F$ is
also a left-adjoint of $F^{\dagger}$, hence $F^{\dagger}$ a right adjoint of $F$.  So if
this prescription gives a left adjoint, it gives a two-sided adjoint.
Next we check that it does.

Suppose $x=(X_i) \in \V^n$ is the 2-vector with vector space $X_i$ in
the $i^{th}$ component, and $y=(Y_j) \in \V^m$ has vector space $Y_j$
in the $j^{th}$ component.  Then $F(x) \in \V^m$ has $j^{th}$ component
$\oplus_{i=1}^n F_{i,j} \otimes X_i$.  Now, a map in $\V^m$ from $F(x)$
to $y$ consists of a linear map in each component, so it is an
$m$-tuple of maps:
\begin{equation}
f_j: \bigoplus_{i=1}^n F_{i,j} \otimes X_i \ra Y_j
\end{equation}
for $j = 1 \dots m$.  But since the direct sum (biproduct) is a
categorical coproduct, this is the same as an $m \times n$ matrix of
maps:
\begin{equation}
f_{ij}: V_{i,j} \otimes X_i \ra Y_j
\end{equation}
for $k = 1 \dots n$ and $j = 1 \dots m$, and $\hom(F(x),y)$ is the
vector space of all such maps.

By the same argument, a map in $\V^n$ from $x$ to $F^{\dagger}(y)$ consists of an
$n \times m$ matrix of maps:
\begin{equation}
g_{ji}: X_i \ra V_{j,i}^{\ast} \otimes Y_j \cong \hom(V_{j,i},Y_j)
\end{equation}
for $i = 1 \dots n$ and $j = 1 \dots m$, and $\hom(x,F^{\dagger}(y))$ is the
vector space of all such maps.

But then we have a natural isomorphism $\hom(F(x),y) \cong \hom(x,F^{\dagger}(y))$
by the duality of $\hom$ and $\otimes$, so in fact $F^{\dagger}$ is a right
adjoint for $F$, and by the above argument, also a left adjoint.

Moreover, no other non-isomorphic matrix defines a 2-linear map with
these properties, and since any functor is naturally isomorphic to
some matrix, this is the sole $F^{\dagger}$ which works.
\end{proof}

\subsection{Example: Group 2-Algebra}

We conclude this section by giving an example of a 2-vector space:

\begin{example}\label{ex:grp2alg} As an example of a KV 2-vector space, consider the
\textit{group 2-algebra} on a finite group $G$, defined by analogy
with the group algebra:

The group algebra $\CG$ consists of the set of elements formed as
formal linear combinations elements of $G$:
\begin{equation}
b = \sum_{g \in G} b_g \cdot g
\end{equation}
where all but finitely many $b_g$ are zero.  We can think of these as
complex functions on $G$.  The algebra multiplication on $\CG$ is
given by the multiplication in $G$:
\begin{equation}
b \star b' = \sum_{g,g' \in G} (b_g b'_{g'}) \cdot g g' \\
\end{equation}
This does not correspond to the multiplication of functions on $G$,
but to \textit{convolution}:
\begin{equation}
(b \star b')_g = \sum_{h\cdot h' = g}b_h b'_{h'}
\end{equation}

Similarly, the \textit{group 2-algebra} $A=\VG$ is the
\textit{category} of $G$-graded vector spaces.  That is, direct
sums of vector spaces associated to elements of $G$:
\begin{equation}
V = \bigoplus_{g \in G} V_g
\end{equation}
where $V_g \in \V$ is a vector space.  This is a $G$-graded vector
space.  We can take direct sums of these pointwise, so that $(V
\oplus V')_g = V_g \oplus V'_g)$, and there is a ``scalar'' product
with elements of $\V$ given by $(W \otimes V)_g = W \otimes V_g$.
There is also a \textit{group 2-algebra} product of $G$-graded vector
spaces, involving a convolution on $G$:
\begin{equation}
(V \star V')_h = \bigoplus_{g \cdot g' = h}V_g \otimes V'_{g'}
\end{equation}

The category of $G$-graded vector spaces is clearly a KV 2-vector
space, since it is equivalent to $\V^k$ where $k = |G|$.  However, it
has the additional structure of a 2-algebra because of the group
operation on the finite set $G$.
\end{example}

\begin{example}
Given a finite group $G$, the category $\catname{Rep(G)}$ has:
\begin{itemize}
\item \textbf{Objects}: Complex representations of $G$ (i.e. functors $\rho : G \ra \V$, where $G$ is seen as a one-object groupoid)
\item \textbf{Morphisms}: Intertwining operators between reps (i.e. natural transformations)
\end{itemize}

This is clearly a 2-vector space generated by the irreducible representations of $G$.
\end{example}

In the next section, we will see that a similar construction shows that the representation categories of finite groupoids are KV 2-vector spaces.  This will be the beginning of our definition of $\FV$.

This highlights one motivation for thinking of 2-vector spaces: the
fact that, in quantum mechanics, one often ``quantizes'' a classical
system by taking the Hilbert space of (square integrable)
$\mathbb{C}$-valued functions on its phase space.  Similarly, one
approach to finding a higher-categorical version of a quantum theory
is to take $\V$-valued functors, as we discuss in more detail in
Section \ref{sec:kv2vsgrpd}.

By restricting our attention to the (essentially) finite case, we
avoid here the analytical issues involved in finding an analog for
$L^2(X)$.

\section{KV 2-Vector Spaces and Finite Groupoids}\label{sec:kv2vsgrpd}

We have now seen that we can get a 2-vector space as a category of
functions from some finite set $S$ into $\V$, and this may have extra
structure if $S$ does.  However, this is somewhat unnatural, since
$\V$ is a category and $S$ a mere set.  It seems more natural to
consider functor categories into $\V$ from some category $\C$.  These
are examples of the generalized 2-vector spaces described by Elgueta
\cite{elgueta}.  Then the above way of looking at a KV 2-vector space
can be reduced to the situation when $\C$ is a discrete category with
a finite set of objects.  However, there are interesting cases where
$\C$ is not of this form, and the result is still a KV vector space.
A relevant class of examples, as we shall show, come from special
kinds of groupoids.

\subsection{Free 2-Vector Space on a Finite Groupoid}

Since we want our 2-vector spaces to have finitely many generators, we
need a condition on the sorts of groupoids we are talking about here.
Of course, since often one works with topological groupoids which may
be uncountable, the kind of finiteness condition we will have to apply
seems restrictive.  A full treatment of, for example, Lie
groupoids, would require much more consideration of infinite
dimensional 2-vector spaces (and indeed 2-Hilbert spaces).  In the
meantime, we can only consider groupoids which are essentially
finite.

We first show that essentially finite groupoids are among the special
categories $\C$ we want to consider:

\begin{lemma}\label{lemma:fgkv}
If $X$ is an essentially finite groupoid, $\catname{Rep(X)} = [\X,\V]$ is a 2-vector space\end{lemma}
\begin{proof}
The groupoid $\X$ is equivalent its skeleton, $\underline{\X}$,
which contains a single object in each isomorphism class.  Since $\X$
is essentially finite, this is a finite set of objects, and each
object has a finite group of endomorphisms.  So
\begin{equation}
X \simeq \coprod_{x \in \underline{\X}} Aut(x)
\end{equation}
where the groups $Aut(x)$ are seen as one-object groupoids.

Then
\begin{align}
[\X,\V] & \simeq \prod_{x \in \underline{\X}} [Aut(x),\V] \\
 & = \prod_{x \in \underline{X}} Rep(Aut(x))
\end{align}

This inherits the biproducts from the categories $Rep(Aut(x))$.  An
irreducible representation of an essentially finite groupoid amounts
to a choice of isomorphism class of objects $[x]$, and an irreducible
representation of the group $Aut(x)$.  By Schur's Lemma, these are
indeed simple objects, since irreducible representations of a group
are simple.

\end{proof}

We notice that we are speaking here of groupoids, and any groupoid
$\X$ is equivalent to its opposite category $\X^{op}$, by an
equivalence that leaves objects intact and replaces each morphism by
its inverse.  So there is no real difference between $[\X,\V]$, the
category of $\V$-valued functors from $\X$, and $[\X^{op},\V]$, the
category of \textit{$\V$-valued presheaves} (or just
``$\V$-presheaves'') on $\X$.  (We also should note that, since our
groupoids are discrete, there is no distinction here between sheaves
and presheaves).

\begin{figure}[h]
\begin{center}
\includegraphics{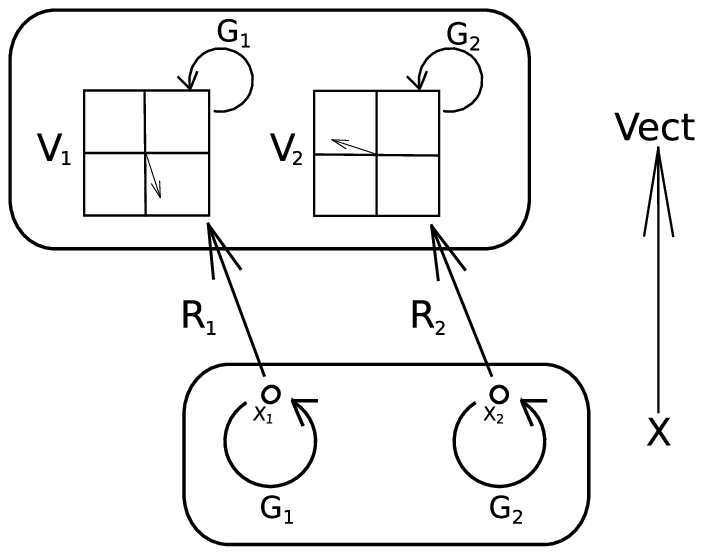}
\end{center}
\caption{\label{fig:vectpresheaf}A $\V$-valued Presheaf on $X$}
\end{figure}

Figure \ref{fig:vectpresheaf} is an illustration of an object in
$[\X,\V]$.

We will use the terminology of ``presheaves'' for objects of $[\X,\V]$
for the sake of highlighting the connection between these results and
the usual facts about presheaves of sets in topos theory - which again
raises questions about topologically interesting groupoids.  This will
be addressed in later work, but for now we consider the algebraic aspect
of the 2-linearization construction by itself.

\subsection{The Ambidextrous Adjunction}\label{sec:ambadj}

Now we want to highlight a result analogous to a standard result for
set-valued presheaves (see, e.g. MacLane and Moerdijk \cite{macmoer},
Theorem 1.9.2).  This is that functors between groupoids induce
2-linear maps between the 2-vector spaces of $\V$-presheaves on them.
For $\Set$-presheaves, there will be a left and a right adjoint to
this functor.  For $\V$-presheaves, these coincide, as we have seen in Theorem \ref{thm:biadjoint} (an inspection of the proof shows that this is essentially because a finite dimensional vector space $V$ is naturally isomorphic to its double dual $V^{\dagger \dagger}$, while the analogous statement is false for sets).  Thus, one says that the
``pushforward'' map is an \textit{ambidextrous adjoint} for the
pullback.  For much more on ambidextrous
adjunctions and their relation to TQFTs, see Lauda \cite{laudaambidjunction}).  This is one important motivation for the present work.  We summarize the statement as follows.

\begin{proposition}\label{thm:2mapadjoints}If $\X$ and $\Y$ are
essentially finite groupoids, a functor $f:\Y \ra \X$ gives two
2-linear maps between KV 2-vector spaces:
\begin{equation}
f^{\ast} : [ \X, \V ] \ra [\Y,\V]
\end{equation}
called ``pullback along $f$'' and
\begin{equation}
f_{\ast} : [\Y,\V] \ra [\X,\V]
\end{equation}
the (two-sided) adjoint to $f^{\ast}$, called ``pushforward along
$f$''
\end{proposition}

\begin{proof}
For any functor $F: \X \ra \V$,
\begin{equation}
f^{\ast}(F) = F \circ f
\end{equation}
which is a functor from $\Y$ to $\V$, the pullback of $F$ along $f$.

To see that this is a 2-linear map, we recall that it is enough to
show it is $\mathbb{C}$-linear, since then biproducts will
automatically be preserved.  But a linear combination of maps in some
$\hom$-category in $[\X,\V]$ is taken by $f^{\ast}$ to the
corresponding linear combination in the $\hom$-category in $[\Y,\V]$,
where maps are now between vector spaces thought of over $y \in \Y$.

So indeed there is a 2-linear map $f^{\ast}$.  But then by Theorem
\ref{thm:biadjoint}, there is a two-sided adjoint of $f^{\ast}$,
denoted $f_{\ast}$.
\end{proof}

\begin{figure}[h]
\begin{center}
\includegraphics{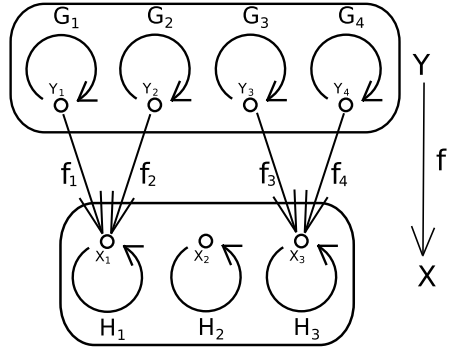}
\end{center}
\caption{\label{fig:grpdfn}A Functor $F : \Y \rightarrow \X$ Between Groupoids}
\end{figure}

In Figure \ref{fig:grpdfn}, we see the essential information contained
in a functor of groupoids.  Any groupoid is equivalent to a skeletal
one (that is, one with just one object in each isomorphism class), so
we illustrate this case.  A skeletal groupoid can be seen as a set of
objects, each labelled by a group.  A functor between groupoids is a
set map, where each ``strand'' of the set map (i.e. each pair
$(y_i,x_j)$ of source and image under the map) is labelled by a
homomorphism $f_i$.  This takes the group $G_i$ of automorphisms of
the source $y_i$ to the group $H_j$ of automorphisms of the target
$x_j$.

It will be useful to have another, more explicit, way to describe the
``pushforward'' map than the matrix-dependent view of Theorem
\ref{thm:biadjoint}.  Fortunately, there is a more intrinsic way to
describe the 2-linear map $f_{\ast}$, the adjoint of $f^{\ast}$.

\begin{definition}\label{def:pushcolimit}
For a given $x \in \X$, the \textit{comma category} $(f \downarrow x)$
has objects which are objects $y \in \Y$ equipped with maps $f(y) \ra
x$ in $\X$, and morphisms which are morphisms $a : y \ra y'$ whose
images make the triangles
\begin{equation}
\xymatrix{
f(y) \ar[d] \ar[r]^{f(a)} & f(y') \ar[dl] \\
x
}
\end{equation}
in $\X$ commute.  Given a $\V$-presheaf $G$ on $\Y$, define
$f_{\ast}(G)(x) = \opname{colim}G(f \downarrow x)$---a colimit in $\V$.

The pushforward of a morphism $b:x \ra x'$ in $\X$, $f_{\ast}(G)(b) :
f_{\ast}(G)(x) \ra f_{\ast}(G)(x')$ is the induced morphism.
\end{definition}

The comma category is the appropriate categorical equivalent of a
\textit{preimage}---rather than requiring $f(y) = x$, one accepts that
they may be isomorphic, in different ways.  So this colimit is a
categorified equivalent of taking a sum over a preimage.  The result
is the \textit{Kan extension} of $G$ along $f$.

Consider the effect of $f_{\ast}$ on a 2-vector $G:\Y \ra \V$ by
describing $f_{\ast}G : \X \ra \V$.  If $F : \X \ra \V$ is as above,
there should be a canonical isomorphism between $[G,f^{\ast}(F)]$ (a
hom-set in $[\Y,\V]$) and $[p_{\ast}(G),F]$ (a hom-set in $[\X,\V]$).

The hom-set $[G,f^{\ast}(F)]$ is found by first taking the pullback of
$F$ along $f$.  This gives a presheaf on $\Y$, namely $F(f(-))$.  The
hom-set is then the set of natural transformations $\alpha : G \ra
f^{\ast}F$.  Given an object $y$ in $\Y$, $\alpha$ picks a
linear map $\alpha_y : F(f(y)) \ra G(y)$ subject to the naturality
condition.

Now, we have seen that, given $f: \Y \rightarrow \X$, this $f_{\ast} :
[\Y,\V] \rightarrow [\X,\V]$ is a 2-linear map, and an ambidextrous
adjoint for $f^{\ast}$.  We would like to describe $f_{\ast}$ more
explicitly.  We shall want to make use of the units and
counits from both the adjunction in which $f_{\ast}$ is a left
adjoint, and that in which it is a right adjoint.  These are described
in the next section.

To describe $f_{\ast}$ in more detail, we use the fact that both $\Y$
and $\X$ are equivalent to unions of finite groups, and so a
$\V$-presheaf on $\Y$ is a functor which assigns a representation of
$Aut(y)$ to each object $y \in \Y$.  Furthermore, if $\Y$ and $\X$ are
skeletal, then $f : \Y \ra \X$ on objects can be any set map, taking
objects in $\Y$ to objects in $\X$.  For morphisms, $f$ gives, for
each object $y \in \Y$, a homomorphism from the group $Aut(y) =
\hom(y,y)$ to the group $Aut(f(y))$.

So the pullback $f^{\ast}$ is fairly straightforward: given $F : \X
\ra \V$, the pullback $f^{\ast}F = F \circ f : \Y \ra \V$ assigns to
each $y \in \Y$ the vector space $F(f(y))$, and gives a representation
of $\opname{Aut}(y)$ on this vector space where $g : y \ra y$ acts by
$f(g)$.  This is the \textit{pullback representation}.  If $f$ is an
inclusion, this is usually called the \textit{restricted
  representation}.  The \textit{pushforward}, or adjoint of pullback,
for an inclusion is generally called finding the \textit{induced
  representation}.  We remark that for the case where $f$ is an
inclusion, Sternberg \cite{sternberg} gives some classical discussion
of this for complex representations, as does Benson \cite{benson} for
more modules over the group ring with a more general base ring $R$.
Here we use the same term for the more general case when $f$ is any
homomorphism.

For any presheaf $F$, the pushforward $f_{\ast} F$ is determined by
the colimit for each component of that essential preimage.  Then for
each $x \in \X$, we first get:
\begin{equation}
\bigoplus_{g : f(y) \ra x} F(y)
\end{equation}
Which is just the direct sum (i.e. biproduct) over the isomorphism
classes in the essential preimage of the corresponding vector spaces.
However, this is not the colimit: an object in the essential preimage
is a pair $(y,g)$, but we note that if $y$ and $y'$ are isomorphic in
$\Y$, such isomorphisms induce isomorphisms of the spaces $F(y)$, and
the colimit will be a quotient which identifies these spaces.  In
general, the colimit will be a direct sum over isomorphism classes
$[y]$ in the essential preimage.  Each term of the sum is isomorphic
to the induced representation of $F(y)$ under the homomorphism
determined by $f$.

Now, consider what the induced representations are for each
isomorphism class.  Any isomorphism class $[y]$ of objects in $\Y$
determines a group $G = Aut(y)$, and similarly $[x] \in
\underline{\X}$ determines $H = Aut(x)$.  So this reduces to the case
where $\Y$ and $\X$ are just groups (seen as one-object categories),
so we have a group homomorphism $f : G \ra H$.  Using the induced
algebra homomorphism $f : \mathbb{C}[G] \ra \mathbb{C}[H]$, one can
directly construct the induced homomorphism as a quotient: $f_{\ast}V
= \mathbb{C}[H] \otimes_{\mathbb{C}[G]} V$.

So for general groupoids, with $V=F(y)$, we have the direct sum:
\begin{equation}\label{eq:inducedrep}
(f_{\ast}F)(x) = \bigoplus_{f(y) \cong x} \mathbb{C}[Aut(x)] \otimes_{\mathbb{C}[Aut(y)]} F(y)
\end{equation}

\begin{figure}[h]
\begin{center}
\includegraphics{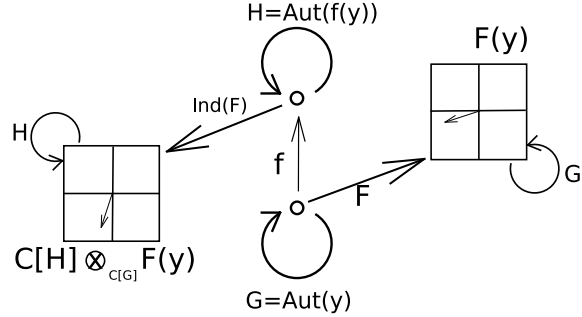}
\end{center}
\caption{\label{fig:indrep}Induced Representation from Homomorphism}
\end{figure}

Figure \ref{fig:indrep} illustrates the induced representation
schematically, for a single object.

\subsection{Units and Counits}\label{sec:units} 

We have observed that the pullback and pushforward maps $f^{\ast}$ and
$f_{\ast}$ are both left and right adjoints.  Thus there are two
adjunctions to consider: $f^{\ast} \dashv f_{\ast}$, where pushforward
is right adjoint to pullback; and $f_{\ast} \dashv f^{\ast}$ where
pushforward is left adjoint to pullback.  For convenience, we refer to
these as the ``right adjunction'' and ``left adjunction''
respectively, after the position of the pushforward.  Each adjunction
has unit and counit, so there are four natural transformations to
describe.  We will identify them as ``right'' and ``left'' unit and
counit following the convention above.  Thus, we have:

\begin{align}
\eta_L : & \Id_{[\Y,\V]} \Longrightarrow f^{\ast} f_{\ast} \\
\epsilon_L : & f_{\ast} f^{\ast} \Longrightarrow \Id_{[\X,\V]} \\
\eta_R : & \Id_{[\X,\V]} \Longrightarrow f_{\ast} f^{\ast} \\
\epsilon_R : & f^{\ast} f_{\ast} \Longrightarrow \Id_{[\Y,\V]}
\end{align}

Once again, it is useful for practical calculations to have a
coordinate-dependent form for these maps, but there is a convenient
intrinsic definition which we shall describe first.  Here again, we
note that Benson \cite{benson} describes the case where $f$ is an
inclusion, in a more general setting than the complex representations
we consider here.

To begin with, we should describe the functors $f^{\ast} f_{\ast}$
(``push-pull''), which is an endofunctor on $[\Y,\V]$, and $f_{\ast}
f^{\ast}$ (``pull-push''), which is an endofunctor on $[\X,\V]$.

For the ``push-pull'', $f^{\ast} f_{\ast}$, we first push a
$\V$-presheaf $F$ on $\Y$ to one on $\X$, then pull back to $\Y$.  On
each object $y \in \Y$, this gives a new presheaf where the vector
space $F(y)$ is replaced by the pullback (i.e. induced representation
of $Aut(y)$) of $f_{\ast}F (f(y))$.  But $f_{\ast}F$ is a presheaf on
$\X$, which, at each $x \in \X$, gives a colimit over the essential
preimage of $x$ in $\Y$, namely $\bigoplus_{[y'] | f(y') \cong x}
\mathbb{C}[Aut(x)] \otimes_{\mathbb{C}[Aut(y')]} F(y')$.  In the case
where $x = f(y)$, this means we get:
\begin{equation}
f^{\ast} f_{\ast} F (x) = \bigoplus_{[y'] | f(y') \cong f(y)} \mathbb{C}[Aut(f(y))] \otimes_{\mathbb{C}[Aut(y')]} F(y')
\end{equation}
thought of as a (left) representation of $Aut(y)$ in the obvious way
(i.e. $g \in Aut(y)$ acts on this space as $f(g)$).

For the ``pull-push'', $f_{\ast} f^{\ast}$, we first pull a
$\V$-presheaf $G$ on $\X$ back to $f^{\ast}G$ on $\Y$.  At each $y \in
\Y$, this assigns the vector space $f^{\ast}G(y) = G(f(y))$ as a
representation of $Aut(y)$. We then push forward to $\X$ to get, at
each $x \in \X$, that:
\begin{equation}
f_{\ast} f^{\ast} G (x) = \bigoplus_{[y] | f(y) \cong x} \mathbb{C}[Aut(x)] \otimes_{\mathbb{C}[Aut(y)]} G(x)
\end{equation}
Note that \textit{a priori} the last space would be $G(f(y))$, but
since $f(y) \cong x$, we have also that $G(f(y)) \cong G(x)$ as
representations of $Aut(y)$.  Here we are implicitly taking a colimit
over the essential preimage of $x$, whose objects are not just $y$
such that $f(y) \cong x$, but rather such $y$ equipped with a
specific isomorphism.  These therefore induce specific isomorphisms of
$G(f(y))$ and $G(x)$, and the quotient implied by the colimit
identifies these spaces.

Now, the description above accords with the usual description of these
functors in the left adjunction.  Since the adjunction is
ambidextrous, it applies in both cases, but to describe the unit and
counit properly, we should note that in general the canonical
description of the left and right adjunctions are different.  (Here
again we note that Benson \cite{benson} shows this for modules over
general rings, which in our case are the group algebras
$\mathbb{C}[Aut(y)]$ etc., in the case of inclusion) We need to take
account of the specific isomorphism between the form we have presented
(natural for the left adjunction), and the form which is natural for
the right adjunction.

The right adjoint is given as:
\begin{equation}
f_{\ast} F (y) = \bigoplus_{[y] | f(y) \cong x} \hom_{\mathbb{C}[Aut(x)]}(\mathbb{C}[Aut(y)],F(y))
\end{equation}
(Note that the case for groups, namely when $\Y$ and $\X$ have only
one object, appears in each term of this direct sum).  The
\textit{Nakayama isomorphism} gives the duality between the two
descriptions of $f_{\ast}$, in terms of $\hom_{\mathbb{C}[Aut(x)]}$
and $\otimes_{\mathbb{C}[Aut(x)]}$, by means of the \textit{exterior
trace map}.  The groupoid case is just the direct sum of group
cases, which looks like:
\begin{equation}
N : \bigoplus_{[y] | f(y) \cong x} \hom_{\mathbb{C}[Aut(y)]}(\mathbb{C}[Aut(x)],F(y)) \ra \bigoplus_{[y] | f(y) \cong x} \mathbb{C}[Aut(x)] \otimes_{\mathbb{C}[Aut(y)]} F(y)
\end{equation}
given by the \textit{exterior trace map} in each factor of the sum:
\begin{equation}\label{eq:nakayama}
N : \bigoplus_{[y] | f(y) \cong x} \phi_y \mapsto \bigoplus_{[y] | f(y) \cong x} \frac{1}{\# Aut(y)}\sum_{g \in Aut(x)} g^{-1} \otimes \phi_y(g)
\end{equation}
Note that the exterior trace map gives an $Aut(x)$-invariant vector,
but the normalization is by the size of $Aut(y)$.  In the case where
the homomorphism is an inclusion, this is interpreted as trace given
by a sum over cosets of $Aut(y)$ in $Aut(x)$, (which is the situation
usually presented in the group case).  We remark here that this factor
will be important in interpreting our 2-functor $\FV$ as a form of
groupoidification.

We can now write down the units and counits explicitly for both
adjunctions in our preferred notation.

The left and right units are natural transformations which, for each
$\V$-presheaf on $\Y$ or $\X$ respectively, gives a morphism which is
itself a natural transformation.  So, in particular the left unit
\begin{equation}\label{eq:leftunit}
\eta_L(F)(y) : F(y) \ra \bigoplus_{[y'] | f(y') \cong f(y)} \mathbb{C}[Aut(f(y))] \otimes_{\mathbb{C}[Aut(y')]} F(y')
\end{equation}
is given by the natural map into the counit:
\begin{equation}
v \mapsto \bigoplus_{[y]} (1 \otimes v)
\end{equation}

Notice the unit map has no contribution in the image from any $y'$
which is not in the isomorphism class $[y]$.  (It is a canonical map
out of the limit which gives the usual form for $f_{\ast}$.)

The right unit map
\begin{equation}\label{eq:rightunit}
\eta_R(G)(x) : G(x) \ra \bigoplus_{[y] | f(y) \cong x} \mathbb{C}[Aut(x)] \otimes_{\mathbb{C}[Aut(y)]} f^{\ast}G(x)
\end{equation}
is found by composing the Nakayama isomorphism (\ref{eq:nakayama})
with the groupoid form of the canonical map for the right adjoint.
This is a direct sum, in which each factor is given by the
multiplication map:
\begin{equation}
v \mapsto \bigl{(} g \mapsto g(v) \bigr{)}
\end{equation}
Thus, the composite is:
\begin{equation}
\eta_R(G)(x) : v \mapsto \bigoplus_{[y] | f(y) \cong x} \frac{1}{\# Aut(y)}\sum_{g \in Aut(x)} g^{-1} \otimes g(v)
\end{equation}

The left and right counits are natural transformations which, for each
$\V$-presheaf on $\X$ or $\Y$ respectively, gives a morphism which is
itself a natural transformation.  So in particular, the left counit
\begin{equation}\label{eq:leftcounit}
\epsilon_L(G)(x) : \bigoplus_{[y] | f(y) \cong x} \mathbb{C}[Aut(x)] \otimes_{\mathbb{C}[Aut(y)]} f^{\ast} G(x) \ra G(x)
\end{equation}
is given by summing multiplication maps:
\begin{equation}
\bigoplus_{[y] | f(y) \cong x} g_y \otimes v \mapsto \sum_{[y] | f(y) \cong x} f(g_y) v
\end{equation}

The right counit map
\begin{equation}\label{eq:rightcounit}
\epsilon_R(F)(y) : \bigoplus_{[y'] | f(y') \cong f(y)} \mathbb{C}[Aut(f(y))] \otimes_{\mathbb{C}[Aut(y')]} F(y') \ra F(y)
\end{equation}
is given by composing the inverse of the Nakayama isomorphism
(\ref{eq:nakayama}) with the evaluation map from the canonical form of
the right adjoint.  Again, the only factor which contributes is $y'
\cong y$, and so we have:
\begin{equation}
\bigoplus_{[y'] | f(y') \cong f(y)}\phi_{y'} \mapsto \phi_y(1)
\end{equation}
So finally, (by using that $\mathbb{C}[Aut(f(y))]$ is canonically isomorphic to its dual using the canonical inner product on the group algebra) the composite is:
\begin{equation}
\epsilon_R(F)(y) : \bigoplus_{[y'] | f(y') \cong f(y)}g_{y'} \otimes v_{y'} \mapsto \frac{\# Aut(y)}{\# Aut(f(y))} g_y(v_y)
\end{equation}
Here we are implicitly using the fact that the objects $y'$ in the
essential preimage come equipped with isomorphisms $f(y') \ra f(y)$
which induce specified isomorphisms $Aut(f(y')) \cong Aut(f(y))$.  In
the colimit which gave the direct sum over isomorphism classes, these
are all naturally identified.

A straightforward check (cancelling the Nakayama isomorphisms) verifies the unit and counit identities:
\begin{align}
(\epsilon_L \cdot Id_{f_{\ast}} ) \circ (Id_{f_{\ast}} \cdot \eta_L) = \Id_{f_{\ast}} \\
(Id_{f^{\ast}} \cdot \epsilon_L) \circ (\eta_L \cdot Id_{f^{\ast}}) = \Id_{f^{\ast}} \\
(\epsilon_R \cdot Id_{f^{\ast}} ) \circ (Id_{f^{\ast}} \cdot \eta_R) = \Id_{f^{\ast}} \\
(Id_{f_{\ast}} \cdot \epsilon_R) \circ (\eta_R \cdot Id_{f_{\ast}}) = \Id_{f_{\ast}}
\end{align}

\section{Spans of Groupoids}\label{sec:spans}

We have already seen how essentially finite groupoids give rise to
2-vector spaces.  In this section, we will show the weak functoriality
of these assignments.  In particular, we first must describe how our
2-functor $\FV$ will produce 2-linear maps from spans of groupoids.

\subsection{2-Linear Maps from Spans of Groupoids}\label{sec:span2lin}

Given a span of groupoids as in Figure \ref{fig:grpd-span}, we can
apply the functor $[ - , \V ]$ to the span diagram
(\ref{eq:spandiag}).  This functor is contravariant, so we get a
cospan:
\begin{equation}\label{eq:zstep3}
\xymatrix{
  & \HV{X} & \\
\HV{A_1} \ar[ur]^{s^{\ast}} & & \HV{A_2} \ar[ul]_{t^{\ast}} \\
}
\end{equation}

We now recall that the pullbacks $s^{\ast}$ and $t^{\ast}$ have
adjoints: this is a direct consequence of Theorem
\ref{thm:2mapadjoints}.  This reveals how to transport a $\V$-presheaf
on $A_1$ along this cospan.  In fact, it gives two 2-linear maps,
which are adjoint.  Thinking of the span as a morphism in
$\Span(\Gpd)$ from $A_1$ to $A_2$, we find a corresponding
2-linear map (though the adjoint is equally well defined).  We first
do a pullback along $s$, giving a $\V$-presheaf on $X$.  Then we use
the adjoint map $t_{\ast}$.  So we have the following:

\begin{definition}\label{def:ZGonS} For a span of groupoids $X : A_1 \ra A_2 \ in \Span (\Gpd)$ define the 2-linear map:
\begin{equation}\label{eq:zstep4}
t_{\ast} \circ s^{\ast} : \HV{A_1} \longrightarrow \HV{A_2}
\end{equation}
\end{definition}

Now, by Theorem \ref{thm:2mapadjoints}, both $s^{\ast}$ and $t_{\ast}$
are 2-linear maps, so the composite $t_{\ast} \circ s^{\ast}$ is also
a 2-linear map.

\begin{remark} We can think of the pullback-pushforward construction
as giving---in the language of quantum field theory---a ``sum over
histories'' for evolving a 2-vector.  Each 2-vector in $\HV{A_1}$ picks
out a vector space for each object of $A_1$.  The 2-linear map we have
described tells us how to evolve this 2-vector along a span.  First we
consider the pullback to $\HV{X}$, which gives us a 2-vector
consisting of all assignments of vector spaces to objects of $X$ which
project to the chosen one in $A_1$.  Each of these objects could be
considered a ``history'' of the 2-vector along the span.  We then
``push forward'' this assignment to $A_2$, which involves a colimit.
This is more general than a sum, though so one could describe this as
a ``colimit of histories''.  It takes into account the symmetries
between individual ``histories'' (i.e. morphisms in $X$).
\end{remark}

So, given a span $X : A_1 \ra A_2$, we can write $\FV(X)$ in terms of its effect on a $\V$-presheaf $G$ on $A_1$, which, at any object $a_2 \in A_2$ gives:
\begin{equation}\label{eq:FVX}
\FV(X)(G)(a_2) = ( \bigoplus_{[x] | t(x) \cong a_2} \mathbb{C}[Aut(a_2)] \otimes_{\mathbb{C}[Aut(x)]} G(s(x)) )
\end{equation}
by exactly the same reasoning as in Section \ref{sec:units}.

It is sometimes useful---particularly when we look at composition of
spans---to break up the direct sum into the contributions from
different objects of $A_1$, like this:
\begin{equation}\label{eq:decomposedFVX}
\FV(X)(G)(a_2) = \bigoplus_{[a_1] \in \underline{A_1}} ( \mathop{\bigoplus_{\stackrel{[x]}{s(x) \cong a_1}}}_{{t(x) \cong a_2}} \mathbb{C}[Aut(a_2)] \otimes_{\mathbb{C}[Aut(x)]} G(a_1) )
\end{equation}

Moreover, there is a convenient way to write down the components of
the 2-linear map associated to a span, which is given by Frobenius
reciprocity.

\begin{proposition}Given basis elements $(a_1,W_1) \in \FV(A_1)$ and
$(a_2,W_2) \in \FV(A_2)$, the matrix elements are:
\begin{equation}\label{eq:frobenius} 
\FV(X)_{(a_1,W_1),(a_2,W_2)} \simeq \bigoplus_{[x]} \opname{hom}_{Rep(\opname{Aut}(x))}[s^{\ast}(W_1),t^{\ast}(W_2)]
\end{equation} \end{proposition}

Here, the direct sum is taken over equivalence classes $[x]$ in the
essential preimage of $(a_1,a_2)$: that is, objects of $X$ mapping to
$a_1$ and $a_2$.  For each $[x]$, the functors $s$ and $t$ define
homomorphisms
\begin{equation}
s_x : \opname{Aut}(x) \ra \opname{Aut}(s(x))
\end{equation}
and
\begin{equation}
t_x : \opname{Aut}(x) \ra \opname{Aut}(t(x))
\end{equation}, which define the induced representations.  We
think of the terms of the direct sum as ``lying over'' the objects
$x$.

So using the adjoint 2-linear map
\begin{equation}
t_{\ast} : \FV(X) \ra \FV(A_2)
\end{equation}
to \textit{push forward} a 2-vector $s^{\ast}F : X \ra \V$ to one
on $A_2$, the above is also, by Frobenius reciprocity:
\begin{equation}
\bigoplus_{[x]} \opname{hom}_{Rep(\opname{Aut}(a_2))}[(t_x)_{\ast} \circ (s_x)^{\ast} W_1 , W_2]
\end{equation}

By Schur's lemma, this says:
\begin{equation}
  \FV(X)(a_1,W_1) = \bigoplus_{[x]} (t_x)_{\ast} \circ (s_x)^{\ast} (a_1,W_1)
\end{equation}
since the components of $\FV(X)(a_1,W_1)$ count the number of copies
of $W_2$ in the pushforward of $W_1$.  (In the remainder of this
paper, we will suppress the subscripts and denote $s_x$ by $s$ and
$t_x$ by $t$ when the context makes clear that we mean the induced
group homomorphism.)

So in fact, $\FV(X)(a_1,W_1)$ is a direct sum of irreducible 2-vectors in
$\FV(A_2)$, given as a sum over $x \in X$ restricting to $a_1, a_2$
of the induced representations along each restriction map.

\subsection{$\FV$ and Composition}

Next we show that $\FV$ preserves horizontal composition of functors
\textit{weakly}---that is, up to a natural isomorphism.  That is, the
composition of the 2-linear maps must be compatible, in a weak sense,
with composition of spans of groupoids.  

To construct the isomorphism explicitly, we look at the weak pullback square in the middle of
(\ref{xy:wkpb}), since the two 2-linear maps being compared
differ only by arrows in this square.  The square as given is a weak
pullback, with the natural isomorphism $\alpha$ ``horizontally''
across the square.  When considering a corresponding square of
categories of $\V$-presheaves, the arrows are reversed.  So, including
the adjoints of $t^{\ast}$ and $T^{\ast}$, namely
$t_{\ast}$ and $T_{\ast}$, we have the square:
\begin{equation}\label{eq:adjunctionmatesquare}
\xymatrix{
   & \HV{X'\circ X} \ar@<1ex>^{T_{\ast}}[rd] & \\
   \HV{X} \ar^{S^{\ast}}[ur] \ar@<1ex>^{t_{\ast}}[rd] & & \HV{X'} \ar@<1ex>^{T^{\ast}}[ul] \\
   & \HV{A_2} \ar_{(s')^{\ast}}[ur] \ar@<1ex>^{t^{\ast}}[ul] &
} \end{equation}

Note that there are two squares here---one by taking only the ``pull''
morphisms $(-)^{\ast}$ from the indicated adjunctions, and the other
by taking only the ``push'' morphisms $(-)_{\ast}$.  The first is just
the square of pullbacks along morphisms from the weak pullback square
of groupoids.  Comparing these is the core of the following
theorem, which gives one of the necessary properties for $\FV$ to be a
weak 2-functor.

We give a more explicit description of the functors $\FV(X' \circ X)$
and $\FV(X') \circ \FV(X)$ below, but remark that this general result
is discussed by Panchadcharam \cite{elango} (Proposition 0.0.1), and
the general theory behind this elaborated on by Street \cite{street2}.

\begin{theorem}\label{thm:FVcomposfunc}The process $\FV$ weakly
preserves composition.  In particular, there is a natural isomorphism
\begin{equation}
\beta_{X',X} : \FV(X' \circ X) \ra \FV(X') \circ \FV(X)
\end{equation}
\end{theorem}
\begin{proof}
Recall that, given the composite of two spans of groupoids in
(\ref{xy:wkpb}), we have 2-linear maps:
\begin{equation}
\FV(X' \circ X) = (t' \circ T)_{\ast} \circ (s \circ S)^{\ast}
\end{equation}
and 
\begin{equation}
\FV(X') \circ \FV(X) = (t')_{\ast} \circ (s')^{\ast} \circ t_{\ast} \circ (s)^{\ast}
\end{equation}
So we want to show there is a natural isomorphism:
\begin{equation}
\beta_{X',X} : (t' \circ T)_{\ast} \circ (s \circ S)^{\ast} \ra
(t')_{\ast} \circ (s')^{\ast} \circ t_{\ast} \circ (s)^{\ast}
\end{equation}
It suffices to show that there is an isomorphism:
\begin{equation}
\gamma: T_{\ast} \circ S^{\ast} \ra (s')^{\ast} \circ t_{\ast}
\end{equation}
between the upper and lower halves of the square in the middle of
(\ref{xy:wkpb}) since then $\beta_{X',X}$ is obtained by tensoring
with identities.

So first taking a $\V$-presheaf $F$ on ${X}$, we get that
$S^{\ast}F$ is a $\V$-presheaf on ${X' \circ X}$.  Now over any
fixed object $x \in X$, we have a set of objects in ${X'
\circ X}$ which restrict to it: there is one for each choice $(g,x')$
which is compatible with $x$ in the sense that $(x,g,x')$ is an object
in the weak pullback - that is, $g: t(x) \ra s'(x')$.  Each
object of this form is assigned
\[ 
S^{\ast}F(x,g,x') = F(x)
\]

Further, there are isomorphisms between such objects, namely pairs
$(h,k)$ as above.  There are thus no isomorphisms except between
objects $(x,g_1,x')$ and $(x,g_2,x')$ for some fixed $x$ and $x'$.
For any such fixed $x$ and $x'$, objects corresponding to $g_1$ and
$g_2$ are isomorphic if
\begin{equation}
g_2 t(h)  = {p'}_1(k) g_1
\end{equation}.
Denote the isomorphism class of any $g$ by $[g]$.

Then we get:
\begin{equation}
T_{\ast} \circ S^{\ast} F (x') = \bigoplus_{x \in X} \Bigl{(}
\bigoplus_{[x']: t(x) \cong s'(x')} \mathbb{C}[Aut({x'})]
\otimes_{\mathbb{C}[Aut(x) \times_{Aut(t(x))} Aut({x'})]} F(x) \Bigr{)}
\end{equation}
since $Aut(x) \times_{Aut(t(x))} Aut({x'})$ is the automorphism group of the object in
${X' \circ X}$ which restricts to $x$ and $x'$.  Notice that although outside direct sum here is written over all
objects $x$ on $S$, the only ones which contribute any factor are
those for which $g : t(x) \ra s'(x')$ for some $g$.  The inside
direct sum is over all isomorphism classes of elements $g$ for which
this occurs: in the colimit, vector spaces over objects with
isomorphisms between them are identified.

Note that in the direct sum over $[g]$, there is a tensor product term
for each class $[g]: t(x) \ra s'(x')$.  By the definition of the tensor
product over an algebra, we can pass elements of $\mathbb{C}[Aut(x)
\times_{Aut(t(x))} Aut({x'})]$ through the tensor product.  These are generated
by pairs $(h,k) \in Aut(x) \times Aut({x'})$ where the images of $h$ and $k$
are conjugate by $g$ so that $t(h) g = g s'(k)$.  These are just
automorphisms of $g$: so this says we are considering objects only up
to these isomorphisms.

This is the result of the ``pull-push'' side of the square applied to
$F$.  Now consider the ``push-pull'' side: $(s')^{\ast} \circ t_{\ast}$.

First, pushing down to $A_2$, we get, on any object $a \in A_2$
\begin{equation}
t_{\ast}F(a) = \bigoplus_{[x] | t(x) \cong a} \mathbb{C}[Aut(a)] \otimes_{\mathbb{C}[Aut(x)]}F(x)
\end{equation}

Then, pulling this back up to $X'$, we find:
\begin{equation}
(s')^{\ast} \circ t_{\ast}F(x') = \bigoplus_{[x] | t(x) \cong s'(x')} \Bigl{(} \mathbb{C}[Aut(t(x))] \otimes_{\mathbb{C}[Aut({x})]} F(x)  \Bigr{)}
\end{equation}
Now we define a natural isomorphism
\begin{equation}
\gamma_{X,X'} :  T_{\ast} \circ S^{\ast} \ra (s')^{\ast} \circ t_{\ast}
\end{equation}
as follows.  For each $x'$, this must be an isomorphism between the above vector
spaces.  The first step is to observe that there is a 1-1
correspondence \textit{between} the terms of the first direct sums, and then
secondly to note that the corresponding terms are isomorphic.

Since the outside direct sums are over all objects $x \in X$ for
which $t(x) \cong s'(x')$, it suffices to get an isomorphism between
each term.  That is, between 
\begin{equation}\label{eq:pullpushcomponentvs}
\mathbb{C}[Aut({x'})] \otimes_{\mathbb{C}[Aut(x) \times_{Aut(t(x))} Aut({x'})]} F(x)
\end{equation}
and
\begin{equation}\label{eq:pushpullcomponentvs}
\mathbb{C}[Aut(t(x))] \otimes_{\mathbb{C}[Aut(x)]} F(x) 
\end{equation}

In order to define this isomorphism, first note that both of these
vector spaces are in fact $\mathbb{C}[Aut({x'})]$-modules.  An element
of $Aut({x'})$ acts on (\ref{eq:pullpushcomponentvs}) in each component
by the standard group algebra multiplication, giving an action of
$\mathbb{C}[Aut({x'})]$ by extending linearly.  An element $l \in
Aut({x'})$ acts on (\ref{eq:pushpullcomponentvs}) by the action of
$s'(l)$ on $\mathbb{C}[Aut(t(x))]$.  Two elements $l_1, l_2 \in [l]$ in the
same equivalence class have the same action on this tensor product,
since they differ precisely by $(h,k) \in Aut(x) \times Aut({x'})$, so that
$l_2 t(h) = s'(k) l_1$.

Also, we notice that, in (\ref{eq:pullpushcomponentvs}), for each $g
\in Aut(t(x))$, the corresponding term of the form $\mathbb{C}[Aut({x'})]
\otimes_{\mathbb{C}[Aut(x) \times_{Aut(t(x))} Aut({x'})]} F(x)$ is generated by
elements of the form $k \otimes v$, for $k \in \mathbb{C}[Aut({x'})]$.  
and $v \in F(x)$.  These are subject to the relations that, for any
$(h,k_1) \in \mathbb{C}[Aut({x})] \times \mathbb{C}[Aut({x'})]$ such that
$t(h) = g^{-1} s'(k_1) g$:
\begin{equation}
kk_1 \otimes v = k (h,k_1) \otimes v = k \otimes (h,k_1) v = k \otimes hv
\end{equation}
since elements of $\mathbb{C}[Aut({x})] \times \mathbb{C}[Aut({x'})]$ act
on $F(x)$ and $\mathbb{C}[Aut({x'})]$  by their projections into
the first and second components respectively.

Now, we define the map $\gamma_{x,x'}$.  First, for any element of the
form $k \otimes v \in \mathbb{C}[Aut({x'})] \otimes_{\mathbb{C}[Aut(x)
\times_{Aut(t(x))} Aut({x'})]} F(x)$ in the $g$ component of the direct sum
(\ref{eq:pullpushcomponentvs}):
\begin{equation}
\gamma_{x,x'} (k\otimes v) = s'(k) g^{-1} \otimes v
\end{equation}
which we claim is in $\mathbb{C}[Aut(t(x))] \otimes_{\mathbb{C}[Aut({x})]}
F(x)$.  This map extends linearly to the whole space.

To check this is well-defined, suppose we have two representatives
$k_1 \otimes v_1$ and $k_2 \otimes v_2$ of the class $k \otimes v$.
So these differ by an element of $\mathbb{C}[Aut(x) \times_{Aut(t(x))} Aut({x'})]$,
say $(h,k)$, so that
\begin{equation}
k_1 = k_2 k
\end{equation}, 
and
\begin{equation}
h v_1 =  v_2
\end{equation}
where
\begin{equation}
t(h) = g s'(k) g^{-1}
\end{equation}

But then
\begin{align}
\gamma_{x,x'}(k_1 \otimes v_1) & = s'(k_1) g^{-1} \otimes v_1 \\
\nonumber & = s'(k_2 k) g^{-1} \otimes v_1 \\
\nonumber & = s'(k_2) g^{-1} g s'(k) g^{-1} \otimes v_1 \\
\nonumber & = s'(k_2) g^{-1} t(h) \otimes v_1\\
\nonumber & = s'(k_2) g^{-1} \otimes h v_1
\end{align}
while on the other hand,
\begin{align}
\gamma_{x,x'} (k_2 \otimes v_2) & = s'(k_2) g^{-1} \otimes v_2\\
\nonumber & = s'(k_2) g^{-1} \otimes h v_1
\end{align}
But these are representatives of the same class in $\mathbb{C}[Aut(t(x))]
\otimes_{\mathbb{C}[Aut({x})]} F(x)$, so $\gamma$ is well defined on
generators, and thus extends linearly to give a well-defined function
on the whole space.

Now, to see that $\gamma$ is invertible, note that given an element $m
\otimes v \in \mathbb{C}[Aut(t(x))] \otimes_{\mathbb{C}[Aut({x})]} F(x)$, we can define
\begin{equation}
\gamma^{-1}(m \otimes v) = 1 \otimes v \in  \bigoplus_{t(x) \ra s'(x')} \mathbb{C}[Aut({x'})]
\otimes_{\mathbb{C}[Aut(x) \times_{Aut(t(x))} Aut({x'})]} F(x)
\end{equation}
in the component coming from the isomorphism class of $g = m^{-1}$ (we
will denote this by $(1 \otimes v)_{m^{-1}}$ to make this explicit,
and in general an element in the class of $g$ will be denoted with
subscript $g$ whenever we need to refer to $g$).

Now we check that this is well-defined.  Given $m_1 \otimes v_1$ and
$m_2 \otimes v_2$ representing the same element $m \otimes v$ of
$\mathbb{C}[Aut(t(x)] \otimes_{\mathbb{C}[Aut({x})]} F(x)$, we must have $h_1
\in Aut(x)$ with
\begin{equation}
m_1 t(h_1) = m_2
\end{equation}
and
\begin{equation}
h_1 v_2 = v_1
\end{equation}
But then applying $\gamma^{-1}$, we get:
\begin{equation}
\gamma^{-1}(m_1 \otimes v_1) = (1 \otimes v_1)_{m_1^{-1}} = (1 \otimes h_1 v_2)_{m_1^{-1}}
\end{equation}
and
\begin{equation}
\gamma^{-1}(m_2 \otimes v_2) = (1 \otimes v_2)_{m_2^{-1}} = (1 \otimes v_2)_{t(h_1)^{-1} m_1^{-1}}
\end{equation}
but these are in the same component, since $g \sim g'$ when $g'
s'(k) = t(h) g$ for some $h \in Aut(x)$ and $k \in Aut({x'})$.  But
then, taking $k=1$ and $h = h_1^{-1}$, we get that $m_1^{-1} \sim
m_2^{-1}$, and hence the component of $\gamma(m \otimes v)$ is well
defined.

But then, consider $m \otimes v = \gamma((k \otimes v)_g) = s'(k)
g^{-1} \otimes v$.  Applying $\gamma^{-1}$ we get:
\begin{equation}
\gamma^{-1} \circ \gamma (k \otimes v)_g = (1 \otimes v)_{g s'(k)^{-1}}
\end{equation}
so we hope that these determine the same element.  But in fact, notice
that the morphism in the weak pullback which gives that $g^{-1}$ and
$s'(k) g^{-1}$ are isomorphic is just labelled by $(h,k) = (1,k)$,
which indeed takes $k$ to $1$ and leaves $v$ intact.  So these are the
corresponding elements under this isomorphism.

So $\gamma$ is invertible, hence an isomorphism.  Thus we define
\begin{equation}
\beta_{X,X'} = 1 \otimes \gamma \otimes 1
\end{equation}
This is the isomorphism we wanted.
\end{proof}

So the $\beta_{X,X'}$ can now be seen as natural transformations explicitly. First consider $\FV({X'}) \circ \FV(X)$, which acts on a presheaf $G$ on $A_1$ as follows.
\begin{align}\label{eq:firststepFV}
\FV(X)(G)(a_2) & = \bigoplus_{[x] | t(x) \cong a_2} \mathbb{C}[Aut(a_2)] \otimes_{\mathbb{C}[Aut(x)]} G(s(x)) ) \\
\nonumber & = \bigoplus_{[a_1] \in \underline{A_1}} (\mathop{\bigoplus_{\stackrel{[x]}{s(x) \cong a_1}}}_{{t(x) \cong a_2}} \mathbb{C}[Aut(a_2)] \otimes_{\mathbb{C}[Aut(x)]} (\FV(X)(G)(a_1))
\end{align}
and then applying $\FV(X')$ to this, we get, rearranging direct sums suitably:
\begin{align}\label{eq:composedFV}
& (\FV({X'}) \circ \FV(X))(G)(a_3) \\
\nonumber = & \bigoplus_{[a_2] \in \underline{A_2}} (\mathop{\bigoplus_{\stackrel{[x']}{s'(x') \cong a_2}}}_{{t'(x') \cong a_3}} \mathbb{C}[Aut(a_3)] \otimes_{\mathbb{C}[Aut({x'})]} (\FV(X)(G)(a_2)) \\
\nonumber = &  \bigoplus_{[a_1] \in \underline{A_1}} \bigl{(} \bigoplus_{[a_2] \in \underline{A_2}} \mathop{\bigoplus_{\stackrel{[x]}{s(x) \cong a_1}}}_{{t(x) \cong a_2}}   \mathop{\bigoplus_{\stackrel{[x']}{s'(x') \cong a_2}}}_{{t'(x') \cong a_3}} \mathbb{C}[Aut(a_3)] \otimes_{\mathbb{C}[Aut({x'})]} \mathbb{C}[Aut(a_2)] \otimes_{\mathbb{C}[Aut(x)]} G(a_1) \bigr{)}
\end{align}
We similarly have:
\begin{equation}\label{eq:FVcomposed}
(\FV(X') \circ \FV(X))(G)(a_3) = \bigoplus_{[a_1] \in \underline{A_1}} ( \mathop{\bigoplus_{\stackrel{[(x,f,x')]}{((s\circ S)(x,f,x') \cong a_1}}}_{{(t'\circ T(x,f,x') \cong a_3}} \mathbb{C}[Aut(a_3)] \otimes_{\mathbb{C}[Aut(x,f,x')]} G(a_1) )
\end{equation}

The isomorphisms $\beta_{X,X'}$ allow us to identify (\ref{eq:composedFV}) and (\ref{eq:FVcomposed}).

\begin{remark}\label{rk:composbeta}
We can also describe the effect of $\beta$ in coordinates - that is, in the
matrix form for a natural transformation of a 2-linear map.  This illustrates the fact that $\mathbb{C}[Aut(a_2)] \cong \bigoplus_{W} W \otimes W^{\star}$, where the sum is over irreducible representations of $Aut(a_2)$.
For suppose
we have a composite of spans, $X' \circ X$.  By Lemma
\ref{lemma:kv2linmatrix}, we have that the functors $T_{\ast} \circ
S^{\ast}$ and $(s')^{\ast} \circ t_{\ast}$ can be written in the form
of a matrix of vector spaces as in (\ref{eq:kv2linmatrix}).

First, $\FV(X' \circ X)$ is given by a matrix indexed by classes of
objects and representations $([a_1],W_1)$ from $A_1$ and $([a_3],W_3)$
from $A_3$.  In the form (\ref{eq:frobenius}), we see that
\begin{align}
& \FV(X' \circ X)_{([a_1],W_1),([a_3],W_3)} \\
\nonumber \simeq & \bigoplus_{[(x,f,x')]} \opname{hom}_{Rep(\opname{Aut}(x,f,x'))}[(s \circ S)^{\ast}(W_1),(t' \circ T)^{\ast}(W_3)]
\end{align}
where $[(x,f,x')]$ represents an equivalence class of objects in the
weak pullback.  

The isomorphisms $\beta_{X,X'}$ take this to the matrix product of
$\FV(X')$ with $\FV(X)$, which has components given by a direct sum
over classes and representations $([a_2],W_2)$ from $A_2$:
\begin{align}\label{eq:cobmatmultcomponent}
& [\FV(X' \circ X)]_{([a_1],W_1),([a_3],W_3)} \\
\nonumber \ralim^{\beta_{X,X'}} & \bigoplus_{([a_2],W)} [\FV(X)]_{([a_1],W_1),([a_2],W_2)} \otimes [\FV(X')]_{([a_2],W_2),([a_3],W_3)}
\end{align}
Recall that $[\FV(X)]_{([a_1],W_1),([a_2],W_2)}$, is a direct sum over
isomorphism classes of objects of $X$ which restrict to $[a_1]$ and
$[a_2]$, with each component being
\begin{equation}
\hom_{\opname{Rep}(Aut(x))}(s^{\ast}W_1,t^{\ast}W_2)
\end{equation}The
component $[\FV(X')]_{([a_2],W_2),([a_3],W_3)}$, is a similar sum
over classes of objects of $X'$ which restrict to $[a_2]$ and $[a_3]$.

The isomorphism $\beta_{X',X}$ identifies the composite, whose
components are sums over objects of $X' \circ X$, with this product.
This $\beta$ consists of isomorphisms in each component.  So in fact,
the $\beta$ are described by their components:
\begin{align}
\label{eq:betaisomorphisms} & \bigoplus_{[(x,f,x')]} \bigl{[} \hom((s\circ S)^{\ast}(W_1),(t'\circ T)^{\ast}(W_3) \bigr{]} \\
\nonumber\stackrel{\beta_{X,X'}}{\longrightarrow} &  \bigoplus_{([a_2],W_2)} \bigl{[} \hom(s^{\ast}(W_1),t^{\ast}(W_2)) \otimes \hom((s')^{\ast}(W_2),(t')^{\ast}(W_3)) \bigr{]} 
\end{align}
Where the second sum is over equivalence classes of $(x,f,x')$ such
that $f : t(x) \ra s'(x')$, and for which $s(x)=a_1$ and
$t'(x')=a_3$.

Since the choice of $[x]$ and $[x']$ amounts to the same thing as the
choice of $[a_2]$, this isomorphism turns a sum over representations
$W_2$ of tensor products of space (of intertwiners), into a sum over
isomorphism (conjugacy) classes of $f \in \opname{Aut}([a_2])$.  This
isomorphism is describing how the representations in the big pullback
decompose.
\end{remark}

\section{Spans of Spans}\label{sec:spanspan}

The situation we are interested in can be represented as an
equivalence class of spans of spans of the following sort:
\begin{equation}\label{eq:spanspan}
\xymatrix{
 & X_1 \ar[dl]_{s_1} \ar[dr]^{t_1} & \\
A_1 & Y \ar[u]_{s} \ar[d]^{t} & A_2  \\
 & X_2  \ar[ul]_{s_2} \ar[ur]^{t_2} & \\
}
\end{equation}

Recall that we assume weak commutativity here - that is, that there
are isomorphisms $\zeta_s : s_1 \circ s \ra s_2 \circ t$ and $\zeta_t :
t_1 \circ s \ra t_2 \circ t$.

Given this situation, which is a 2-morphism for the associated bicategory
of spans, we want to get a 2-morphism in the bicategory $\iiV$.  That
is to say, a natural transformation $\alpha=\FV(Y)$ between a pair of
2-linear maps.  In this section, we show how to construct $\FV(Y)$.

\subsection{2-Morphisms from Spans of Spans}

We begin by noting that the diagram (\ref{eq:spanspan}), which weakly
commutes up to isomorphisms $\zeta_s$ and $\zeta_t$, gives rise to a
diagram of pullback functors: 
\begin{equation}\label{eq:spanspanpullback}
\xymatrix{
 & \HV{X_1} \ar[d]_{s^{\ast}} & \\
\HV{A_1} \ar[dr]_{{s_2}^{\ast}} \ar[ur]_{{s_1}^{\ast}} & \HV{Y} & \HV{A_2} \ar[ul]^{{t_1}^{\ast}} \ar[dl]^{{t_2}^{\ast}} \\
 & \HV{X_2} \ar[u]^{t^{\ast}} & \\
}
\end{equation}
which commutes up to isomorphisms:
\begin{equation}
\zeta_s^{\ast} : s^{\ast} \circ s_1^{\ast}(V) \tilde{\ra} t^{\ast} \circ s_2^{\ast}(V)
\end{equation}
and similarly 
\begin{equation}
\zeta_t^{\ast} : s^{\ast} \circ t_1^{\ast}(V) \tilde{\ra} t^{\ast} \circ t_2^{\ast}(V)
\end{equation}

We want to get a natural transformation from $\FV(X_1)$ and $\FV(X_2)$
from this diagram.  In section \ref{sec:coordY} we show how this can
be described as a ``pull-push'' process, similar to the one used to
define the 2-linear maps, but first it can be defined in terms of the
unit and counit maps we have already defined.

\begin{definition}\label{def:Von2Span}Given a span between spans,
$Y : X_1 \ra X_2$, for $X_1,X_2 : A_1 \ra A_2$, then
\begin{equation}
\FV(Y) : \FV(X_1) \ra \FV(X_2)
\end{equation} is the natural transformation given as
\begin{equation}
\FV(Y) = \epsilon_{L,t} \circ ((\zeta_t)_{\ast} \otimes ((\zeta_s)^{\ast})^{-1} ) \circ \eta_{R,s} : (t_1)_{\ast} s_1^{\ast} \Longrightarrow (t_2)_{\ast} s_2^{\ast}
\end{equation}
where $\epsilon_{L,t}$ is the counit (\ref{eq:leftcounit}) for the
left adjunction associated to $t$, and $\eta_{R,s}$ is the unit
(\ref{eq:rightunit}) for the right adjunction associated to $s$.
\end{definition}

We comment here that this composition of left unit followed by right
counit can be interpreted as a ``pull-push'' in a sense that can be
seen more precisely when we consider this natural transformation in
coordinates.  In the special case where $A_i = 1$, this recovers
groupoidification in the sense of Baez and Dolan, as shown in Theorem
\ref{thm:groupoidification}.

\begin{remark}Henceforth, we will assume that the diagram (\ref{eq:spanspan})
commutes strictly - that is, $\zeta_s$ and $\zeta_t$ are identity
2-morphisms.  This is a mild assumption, since the diagram can always
be ``strictified'' by taking equivalent groupoids for source and
target which are skeletal.  Similar arguments will follow through if
not, but in this simpler case, we simply have:
\begin{equation}
\FV(Y) = \epsilon_{L,t} \circ \eta_{R,s} : (t_1)_{\ast} s_1^{\ast} \Longrightarrow (t_2)_{\ast} s_2^{\ast}
\end{equation}
\end{remark}

For our construction to give a 2-functor, this must agree with
composition in two ways. The first is strict preservation of vertical
composition; the second is preservation of horizontal composition as
strictly as possible (i.e. up to the isomorphisms $\beta$ which make
comparison possible - as we will see).  We will show in 

\begin{lemma}\label{thm:vertcompos}The assignment $\FV(Y)$ to spans of spans
given in Definition \ref{def:Von2Span} preserves vertical composition strictly:
$\FV(Y' \circ Y) = \FV(Y') \circ \FV(Y)$.
\end{lemma}
\begin{proof}
Suppose we have a vertical composite of two spans between spans, here
written as 2-cells:
\begin{equation}\label{xy:verticalcomposition}
\xymatrix{
A_1 \ar@/^3pc/[rrr]^{X_1}="1" \ar[rrr]^{X_2}="2" \ar@/_3pc/[rrr]_{X_3}="3" & & & A_2 \\
 {\ar@{=>}"1"+<0ex,-2.5ex> ;"2"+<0ex,+2.5ex>^{Y}}
 {\ar@{=>}"2"+<0ex,-2.5ex> ;"3"+<0ex,+2.5ex>^{Y'}}
}
\end{equation}

The composition is given by a weak pullback (taken up to isomorphism)
- that is, a diagram of the form (\ref{xy:wkpb}), with the $Y$ and
$Y'$ in place of $X$ and $X'$, and the $X_i$ in place of the $A_i$.
We use the same notation for the maps in all these spans.  Of course,
each object in (\ref{xy:verticalcomposition}) comes equipped with
(commuting) maps into $A_1$ and $A_2$, but we can ignore these here.

So the source and target maps for $Y$ are $s$ and $t$, and those for
$Y'$ are $s'$ and $t'$, and we can easily write:
\begin{equation}
\FV(Y') \circ \FV(Y) = \epsilon_{L,t'} \circ \eta_{R,s'} \circ \epsilon_{L,t} \circ \eta_{R,s}
\end{equation}

Now, to write $\FV(Y' \circ Y)$, we recall that the vertical composite
is formed by weak pullback of spans, with the resulting source and
target maps $s \circ S$ and $t' \circ T$, where the groupoid in the
span $Y' \circ Y$ is the comma category whose objects are of the form
$(y',g_2,y)$, with $g_2 : t(y) \ra s'(y')$ in $X_2$, and $S$ and $T$
are the natural projections onto $Y$ and $Y'$.  Then of course
\begin{equation}
\FV(Y' \circ Y) = \epsilon_{L,(t' \circ T)} \circ \eta_{R,(s \circ S)}
\end{equation}

Now, a composite of adjunctions is an adjunction (see for instance
MacLane \cite{maclane} IV.8), and the unit and counit of the composite is
given in a standard way, so we have:
\begin{equation}
\epsilon_{L,(t' \circ T)} = \epsilon_{L,t'} \circ (\Id_{(t')_{\ast}} \otimes \epsilon_{L,T} \otimes \Id_{(t')^{\ast}})
\end{equation}
and
\begin{equation}
\eta_{R,(s \circ S)} = (\Id_{s_{\ast}} \otimes \eta_{R,S} \otimes \Id_{s^{\ast}}) \circ \eta_{R,s}
\end{equation}

So we get
\begin{equation}
\FV(Y' \circ Y) = \epsilon_{L,t'} \circ (\Id_{(t')_{\ast}} \otimes \epsilon_{L,T} \otimes \Id_{(t')^{\ast}}) \circ (\Id_{s_{\ast}} \otimes \eta_{R,S} \otimes \Id_{s^{\ast}}) \circ \eta_{R,s}
\end{equation}

So to get strict composition, we just need that
\begin{equation}
(\Id_{(t')_{\ast}} \otimes \epsilon_{L,T} \otimes \Id_{(t')^{\ast}}) \circ (\Id_{s_{\ast}} \otimes \eta_{R,S} \otimes \Id_{s^{\ast}}) = \eta_{R,s'} \circ \epsilon_{L,t}
\end{equation}
This follows the same pattern as the proof for the fact that $\FV$
weakly preserves composition of morphisms.  Note that the two sides of
this expression are the top and bottom of a (weak) pullback square.
So in particular, the argument for weak preservation of composition of
spans shows that we also have a pullback square for the induced
functors.

In particular (ignoring the identity maps), we first get the right
unit for $S : Y' \circ Y \ra Y$, which at $y \in G(y) \cong G(a_1)$ gives:
\begin{equation}
\eta_{R,S} : v \mapsto \bigoplus_{[(y,f,y')]} \frac{1}{\# Aut(y,f,y')} \sum_{g \in Aut(y)} (g^{-1}) \otimes g(v)
\end{equation}
and the left counit for $T : Y' \circ Y \ra Y'$ takes this to.  So this gives
\begin{equation}
\epsilon_{L,T} \circ \eta_{R,S}: v \mapsto \sum_{[(y,f,y')]} \frac{1}{\# Aut(y,f,y')} \sum_{g \in Aut(y)} (g^{-1}) \otimes g(v)
\end{equation}
since the target space is now already a tensor product over $\mathbb{C}[Aut(x_2)]$.  Similarly, on the other side we have first the left counit
for $t : Y \ra X_2$, then the right unit for $s' : Y' \ra X_2$, giving:
\begin{equation}
\eta_{R,s'} \circ \epsilon_{L,t} : v \mapsto \sum_{[y'] | s'(y') \cong t(y)} \frac{1}{\# Aut(t(y))} \sum_{g \in Aut(x_2)} g^{-1} \otimes g(v)
\end{equation}
Since the sources and targets of the maps (\ref{eq:leftcounit}) and
(\ref{eq:rightunit}) are in a pullback square, the coefficients
arising from the Nakayama isomorphisms will yield the same group averages, and the terms
of the implied direct sum correspond pairwise.  So these maps are indeed equal.
\end{proof}

We also must show that $\FV$ respects horizontal composition weakly.  To make this clear, it will be convenient to write source and target 2-linear maps in the form (\ref{eq:composedFV}) and (\ref{eq:FVcomposed}).

\begin{lemma}\label{thm:horizcompos}The assignment $\FV(Y)$ to spans of spans given by
Definition \ref{def:Von2Span} preserves horizontal composition strictly, up to
the isomorphism weakly preserving composition of the source and target
morphisms:
\begin{equation}\label{eq:horizcompos}
\xymatrix{
\FV(A_1) \ar@/^3pc/[rrrr]^{\FV({X'}_1 \circ X_1)}="a" \ar@{}[rrrr]^{}="c" \ar@/_3pc/[rrrr]_{\FV({X'}_2 \circ X_2)}="b" \ar@/^1pc/[rr]^{\FV(X_1)}="1" \ar@/_1pc/[rr]_{\FV(X_2)}="2"  & & \FV(A_2') \ar@/^1pc/[rr]^{\FV(X'_1)}="3" \ar@/_1pc/[rr]_{\FV({X'}_2)}="4" & & \FV(A_3) \\
 {\ar@{=>}"1"+<0ex,-2.5ex> ;"2"+<0ex,+2.5ex>^{\FV(Y)}}
 {\ar@{=>}"3"+<0ex,-2.5ex> ;"4"+<0ex,+2.5ex>^{\FV(Y')}}
 {\ar@{=>}"a"+<0ex,-2.5ex> ;"c"+<0ex,+2.5ex>^{\beta_{X_1,{X'}_1}}}
 {\ar@{=>}"c"+<0ex,-2.5ex> ;"b"+<0ex,+2.5ex>^{\beta_{X_2,{X'}_2}^{-1}}}
}
=
\xymatrix{
\FV(A_1) \ar@/^2pc/[rrrr]^{\FV({X'}_1 \circ X_1)}="1" \ar@/_2pc/[rrrr]_{\FV({X'}_2 \circ X_2)}="2" & & & & \FV(A_3) \\
 {\ar@{=>}"1"+<0ex,-2.5ex> ;"2"+<0ex,+2.5ex>^{\FV(Y' \circ Y)}}
}
\end{equation}
\end{lemma}
\begin{proof}
(Elsewhere, we have used the same notation for horizontal and vertical
  composition of all kinds, to simplify notation and because context
  made this unambiguous.  In this proof it will be helpful to
  distinguish the two, so we write vertical composition with no
  symbol, concatenating natural transformations between 2-linear
  maps.)

Begin by writing the spans explicitly.  The situation for a
horizontal composite of 2-morphisms in $\Span(\Gpd)$ looks like:
\begin{equation}\label{xy:horizontalcomposite}
\xymatrix{
 & & {X'}_1 \circ X_1 \ar[ld]_{S_1} \ar[rd]^{T_1} & & \\
 & X_1 \ar@{=>}[rr]^{\alpha_1}_{\sim} \ar[ld]_{s_1} \ar[rd]^{t_1} & & X'_1 \ar[ld]_{s'_1} \ar[rd]^{t'_1} & \\
A_1 & Y \ar[u]^{s} \ar[d]_{t} \ar@{-->}[l]_{\sigma} \ar@{-->}[r]^{\tau} & A_2 & Y' \ar[u]^{s'} \ar[d]_{t'} \ar@{-->}[l]_{\sigma '} \ar@{-->}[r]^{\tau '} & A_3 \\
 & X_2 \ar@{=>}[rr]^{\alpha_2}_{\sim} \ar[lu]^{s_2} \ar[ur]_{t_2} & & X'_2 \ar[lu]^{s'_2} \ar[ur]_{t'_2} & \\
 & & {X'}_2 \circ X_2 \ar[lu]^{S_2} \ar[ru]_{T_2} & &
}
\end{equation}
(Note that here again we are assuming the 2-morphisms $Y$ and $Y'$ are
strict, so that $s_1 \circ s = s_2 \circ t$ and similarly for the
other composites.  We represent these by the dotted arrows $\sigma$,
$\tau$, $\sigma '$ and $\tau '$.  As before, a similar argument would
go through if these spans of span maps were only weakly commuting, but
we would need the $\zeta$ natural transformations as discussed for
(\ref{xy:spanmapspan})).

Now, the functor $\FV$ assigns 2-linear maps to the spans $X_1$,
$X_2$, ${X'}_1$, and ${X'}_2$, and their composites, and natural
transformations to $Y$ and $Y'$.  Then the horizontal composite is a
natural transformation between 2-linear maps:
\begin{equation}\label{eq:horizFV}
\FV(Y') \circ \FV(Y) : \FV({X'}_1) \circ \FV(X_1) \ra \FV({X'}_2) \circ \FV(X_2)
\end{equation}
And we can calculate as in the proof of Theorem \ref{thm:vertcompos} that:
\begin{equation}
 \FV(Y')(a_1) = \epsilon_{L,t'} \eta_{R,s'} : v \mapsto \bigoplus_{[y] | s'(y')  \cong x_1}\frac{\# Aut(y)}{\# Aut(y,f,y')} v
\end{equation}
and
\begin{equation}
 \FV(Y) = \epsilon_{L,t} \eta_{R,s}
\end{equation}

So the composite is just
\begin{equation}
\FV(Y') \circ \FV(Y) = (\epsilon_{L,t'} \eta_{R,s'})\circ (\epsilon_{L,t} \eta_{R,s})
\end{equation}

We recall that $\FV(X') \circ \FV(X)$ is described explicitly in
(\ref{eq:composedFV}) and $\FV(X' \circ X)$ in (\ref{eq:FVcomposed}).
Finding these for the $X_i$ gives a total of four functors here.  We
next describe natural transformations between these.

As shown in Theorem \ref{thm:FVcomposfunc}, there are comparison
isomorphisms
\begin{equation}
\beta_{X_i,{X'}_i} : \FV({X'}_i) \circ \FV(X_i) \ra \FV({X'}_i \circ X_i)
\end{equation}
which will necessarily be involved in the isomorphism we are looking
for.  These derive from the $\alpha$ isomorphisms in the weak pullback
in ${X'}_i \circ X_i$.

Composing with these comparison isomorphisms as in
(\ref{eq:horizcompos}) gives:
\begin{equation}
(\Id_{s_1^{\ast}} \circ \beta_{X_1,{X'}_1} \circ \Id_{({t'}_1)^{\ast}}) (\epsilon_{L,t'} \eta_{R,s'})  (\epsilon_{L,t} \eta_{R,s}) (\Id_{{s_2}^{\ast}} \circ (\beta_{X_2,{X'}_2})^{-1} \circ \Id_{({t'}_2)^{\ast}})
\end{equation}

Now, the $\beta$ isomorphisms simply allow us to identify the spaces here, so it suffices to describe the maps, and in particular the coefficients which arise.  At any presheaf $G$ on $A_1$, in the summand for $[a_1] \in \underline{A_1}$ and $[x'_1] \in X'_1$:
\begin{equation}
 (\epsilon_{L,t} \eta_{R,s}) : v \mapsto \sum_{[y] | s(y) \cong x_1} \frac{1}{\# Aut(y)} \sum_{g \in Aut(x_1)} g^{-1} \otimes g(v)
\end{equation}
which is a map between spaces of the form (\ref{eq:firststepFV}) associated to $X_1$ and $X_2$.  Now this becomes the $v'$ when we take the full map between spaces like (\ref{eq:composedFV}), where we have:
\begin{equation}
  (\epsilon_{L,t'} \eta_{R,s'}): v' \mapsto \sum_{[y'] | s'(y) \cong {x'}_1} \frac{1}{\# Aut(y')} \sum_{h \in Aut({x'}_1)} h^{-1} \otimes h(v')
\end{equation}
so finally we get:
\begin{equation}
v \mapsto \sum_{[y'] | s'(y) \cong {x'}_1} \frac{1}{\# Aut(y')} \sum_{h \in Aut({x'}_1)} h^{-1} \otimes h \bigl{(} \sum_{[y] | s(y) \cong x_1} \frac{1}{\# Aut(y)} \sum_{g \in Aut(x_1)} g^{-1} \otimes g(v) \bigr{)}
\end{equation}
which we want to show is the same as the natural transformation
associated to $\FV(Y' \circ Y)$:
\begin{equation}\label{eq:FVhoriz}
\FV(Y' \circ Y) = \epsilon_{L,(t,t')} \circ \eta_{R,(s,s')} : \FV({X'}_1 \circ X_1) \ra \FV({X'}_2 \circ X_2)
\end{equation}
which is a map between two spaces of the form (\ref{eq:FVcomposed}).

This $Y' \circ Y : {X'}_1 \circ X_1 \ra {X'}_2 \circ X_2$ is a span of
span maps which is given as follows.  We take the horizontal composite
of the spans $A_1 \stackrel{\sigma}{\la} Y \stackrel{\tau}{\ra} A_2$
and $A_2 \stackrel{\sigma '}{\la} Y' \stackrel{\tau '}{\ra} A_3$.
This is a weak pullback taken up to isomorphism.  The pullback square
commutes weakly, say up to $\xi$.  Then the groupoid $Y' \circ Y$ has
maps into $Y$ and $Y'$, and therefore by composition with $s$ and
$s'$, it has maps into $X_1$ and ${X'}_1$.  By the universal property
of the weak pullback ${X'}_1 \circ X_1$, there is a map $S : Y' \circ
Y \ra {X'}_1 \circ X_1$.  Similarly, there is $T : Y' \circ Y \ra
{X'}_2 \circ X_2$.

We can see what this is by taking the weak pullback giving $Y' \circ
Y$, which we take to be the comma category whose objects are of the
form $(y,f,y')$ where $f : \tau(y) \ra \sigma '(y')$ in $A_2$.  Then
the $S$ and $T$ given by the universal property are just $S = (s,s')$
giving objects like $(x_1,f,{x'}_1)$ and $T = (t,t')$ giving objects
like $(x_2,f,{x'}_2)$.  In particular, the morphism $f \in A_2$ is
left intact.  (Different isomorphic choices for weak pullback could of
course change $f$).

Given this, we have, in the summand for a given $[a_1] \in \underline{A_1}$:
and a particular $[(x_1,f,x'_1)] \in X'_1 \circ X_1$, we have:
\begin{equation}
\eta_{R,(s,s')} : v \mapsto \mathop{\bigoplus_{\stackrel{[y,f,y']}{s(y) \cong x_1}}}_{s'(y') \cong x'_1} \frac{1}{\# Aut(y,f,y')} \sum_{(g,h) \in Aut(x,f,x')} (g,h)^{-1} \otimes (g,h)(v)
\end{equation}
and then
\begin{equation}
\epsilon_{L,(t,t')} \eta_{R,(s,s')} : v \mapsto \mathop{\sum_{\stackrel{[y,f,y']}{s(y) \cong x_1}}}_{s'(y') \cong x'_1}  \frac{1}{\# Aut(y,f,y')} \sum_{(g,g') \in Aut(x,f,x')} (g,h)^{-1} \otimes (g,h)(v)
\end{equation}

So in fact, since we are in a weak pullback square, the size of the automorphism groups in the two expressions we have found will be in the same ratios, and so it becomes clear that, using the $\beta$ isomorphisms as seen in the proof of Theorem \ref{thm:FVcomposfunc}:
\begin{equation}
\FV(Y' \circ Y) = \beta_{X_1,{X'}_1} (\FV(Y') \FV(Y)) (\beta_{X_2,{X'}_2})^{-1}
\end{equation}
as required.
\end{proof}

\subsection{Coordinate Description of 2-Morphisms}\label{sec:coordY}

In this section, we discuss the behaviour of $\FV$ on 2-morphisms,
namely the assignment of a natural transformation to a span of span
maps.  As discussed in Section \ref{sec:kv2vsgrpd}, any natural
transformation between a pair of 2-linear maps between KV 2-vector
spaces can be represented as a matrix of linear operators, as in
(\ref{eq:kvnattransmatrix}).  We would like to describe explicitly the
linear maps composing $\FV(Y)$ and some consequences.

To motivate the rest, we can begin with the special case of
$\hom_{\Span(\Gpd)}(\catname{1},\catname{1})$, where $\catname{1}$ is
the trivial groupoid with one object (which we denote $\star$) and its
identity morphism.  We can summarize the effect of $\FV$ on this
$\hom$-category by the following theorem:

\begin{theorem}\label{thm:groupoidification}
On $\hom_{\Span(\Gpd)}(\catname{1},\catname{1})$, the 2-functor $\FV$,
expressed in coordinates, reproduces groupoidification in the sense of
(\ref{eq:baezgpdify}).
\end{theorem}
\begin{proof}
First, we note that $\FV(\catname{1}) \cong \V$, since the only
irreducible representation of the trivial group (basis object) is
$\mathbb{C}$ itself.

Since $\catname{1}$ is terminal in $\Gpd$, any groupoid has a unique
map into it.  Thus, any span from $\catname{1}$ to $\catname{1}$ is of
the form:
\begin{equation}
\catname{1} \stackrel{!}{\la} X \stackrel{!}{\ra} \catname{1}
\end{equation}
which just amounts to a choice of $X$.  Then $\FV(X) = !_{\ast}
!^{\ast}$ can be described as a $1 \times 1$ matrix of vector spaces,
\begin{align}
\FV(X)_{(\star,\mathbb{C}),(\star,\mathbb{C})} & = \hom(!^{\ast}\mathbb{C},!^{\ast}\mathbb{C}) \\
 & \cong \bigoplus_{[x] \in \underline{X}} \mathbb{C}
\end{align}
since $!^{\ast}\mathbb{C}$ is the representation of $X$ assigning a
copy of $\mathbb{C}$ to each object.  In particular, inducing up the
representation $\mathbb{C}$ gives, at each $x \in X$, the representation
\begin{equation}
\mathbb{C}[Aut(\star)] \otimes_{\mathbb{C}[Aut(x)]} \mathbb{C} = \mathbb{C}
\end{equation}
since $!^{\ast}\mathbb{C}$, is the trivial representation of $Aut(x)$.

For each isomorphism class in
$\underline{X}$, we thus get a copy of $\hom{\mathbb{C},\mathbb{C}
  \cong \mathbb{C}}$.  This is the vector space associated to $X$ by
  groupoidification.

Similarly, a 2-morphism $Y : X_1 \ra X_2$ just amounts to an
isomorphism class of spans of groupoids (since $Y$ and the $X_i$ have
unique maps to $\catname{!}$.  Then the linear map
\begin{equation}
\FV(Y)_{(\star,\mathbb{C}),(\star,\mathbb{C})} : \FV(X_1)_{(\star,\mathbb{C}),(\star,\mathbb{C})} \ra \FV(X_2)_{(\star,\mathbb{C}),(\star,\mathbb{C})}
\end{equation}
just becomes a map
\begin{equation}
T(Y) : \mathbb{C}[\underline{X_1}] \ra \mathbb{C}[\underline{X_2}]
\end{equation}
given by $T(Y) = \epsilon_{L,t} \circ \eta_{R,s}$.  

By the above, (\ref{eq:rightunit}), using $F(x_1) \cong \mathbb{C}$, can
be written:

\begin{equation}
\eta_{R,s}(F)(y) : \mathbb{C} \ra \bigoplus_{[y] | s(y) \cong x_1} \mathbb{C}[Aut(x_1)] \otimes_{\mathbb{C}[Aut(y)]} \mathbb{C}
\end{equation}
So now for each $[x_1] \in \underline{X_1}$, every $y$ in the
essential preimage of $x_1$ under $s$ gets a copy of the trivial
representation $\mathbb{C}$ for each coset of $Im(Aut(y))$ in
$Aut((x_1)$.  This describes a decomposition of a representation of
$Aut(y)$ in terms of irreps (all of which are necessarily trivial in
this case).  Call this representation $G$.  In particular, a vector in
$\FV{X_1}_{(\star,\mathbb{C}),(\star,\mathbb{C})}$ gives a complex
number at each $[x_1]$.  The unit $\eta_{R,s}$ takes such a vector $v$
to, at each $y$ with $s(y) \cong x_1$,
\begin{equation}
\frac{1}{\# Aut(y)} \sum_{g \in Aut(x_1)} g^{-1} \otimes 1
\end{equation}

By commutativity for the span of span maps (which is necessarily
strict here!), we also must have that
\begin{equation}
\bigoplus_{[y'] | t(y') \cong t(y)} \mathbb{C}[Aut(t(y))] \otimes_{\mathbb{C}[Aut(y')]}G(y') \cong \bigoplus_{[y'] | s(y') \cong s(y)} \mathbb{C}[Aut(s(y))] \otimes_{\mathbb{C}[Aut(y)]}\mathbb{C}
\end{equation}

Similarly, then, using this  (\ref{eq:rightcounit}) can be written:
\begin{equation}
\epsilon_{L,t}(G)(y) : \bigoplus_{[y'] | t(y') \cong t(y)} \mathbb{C}[Aut(t(y))] \otimes_{\mathbb{C}[Aut(y')]}G(y') \ra G(y)
\end{equation}

So now consider the vector $v \in
\FV(X_1)_{(\star,\mathbb{C}),(\star,\mathbb{C})}$ which gives $1$ at
$[x_1]$ and $0$ elsewhere.  (That is, it gives the identity
intertwining map between the copies of the representation
$!^{\ast}\mathbb{C}$ at objects in $[x_1]$ and the zero intertwiner
elsewhere).  Then the natural transformation induces a map on the
coefficient:
\begin{equation}
\eta_{R,s}: v \mapsto \bigoplus_{[y] | s(y) \cong x_1} \frac{1}{\# Aut(y)} \sum_{g \in Aut(x_1)} g^{-1} \otimes 1
\end{equation}
but then suppose we look for the coefficient of the result at $[x_2]
\in \underline{X_2}$.  Only those $y$ over $[x_2]$ will contribute, but
then, since the $g^{-1}$ have no effect on vectors in $\mathbb{C}$, we get:
\begin{equation}
\epsilon_{L,t} : \bigoplus_{[y] | s(y) \cong x_1} \frac{1}{\# Aut(x_1)} \sum_{g \in Aut(x_1)} g^{-1} \otimes 1 \mapsto \sum_{y | (s,t)(y) \cong (x_1,x_2)} \frac{\# Aut(x_1)}{\# Aut(y)}
\end{equation}

But this is just
\begin{equation}
\# Aut(x_1) \sum_{y | (s,t)(y) \cong (x_1,x_2)} \frac{1}{\# Aut(y)} = \# Aut(x_1) | \widehat{(x_1,x_2)} |
\end{equation}
where the second term is the groupoid cardinality of the essential
preimage of $(x_1,x_2)$.  This is just the coefficient we find in
groupoidification in the sense of Baez and Dolan.
\end{proof}

Similar calculations apply for less trivial situations as well,
although for these we will require a little more of the representation
theory of the groupoids $A_i$.

\begin{lemma}\label{lemma:coordFVY}Given a (strict) span between spans,
$Y : X_1 \ra X_2$, for $X_1,X_2 : A_1 \ra A_2$, then the natural transformation
\begin{equation}
\FV(Y) : \FV(X_1) \ra \FV(X_2)
\end{equation} is a natural transformation given by a matrix of linear
operators:
\begin{equation}
\FV(Y)_{([a_1],W_1),([a_2],W_2)} : \FV(X_2)_{([a_1],W_1),([a_2],W_2)} \ra \FV(X_2)_{([a_1],W_1),([a_2],W_2)} 
\end{equation}
or equivalently
\begin{align}\label{eq:FVYhomsets}
\FV(Y)_{([a_1],W_1),([a_2],W_2)} & : \bigoplus_{[x_1]}\opname{hom}_{Rep(\opname{Aut}(x_1))}[s_1^{\ast}(W_1),t_1^{\ast}(W_2)]\\ 
\nonumber & \ra \bigoplus_{[x_2]}\opname{hom}_{Rep(\opname{Aut}(x_2))}[s_2^{\ast}(W_1),t_2^{\ast}(W_2)]
\end{align}

Such that for each block $([x_1],[x_2])$, the corresponding
linear operator behaves as follows: for $f \in
\opname{hom}[s_1^{\ast}(W_1),t_1^{\ast}(W_2)]$ we get:
\begin{equation}
\FV(Y)_{([a_1],W_1),([a_2],W_2)}|_{(x_1,x_2)}(f) = |\widehat{(x_1,x_2)}| \sum_{g \in \opname{Aut}(x_1)} g^{-1} f g
\end{equation}
where $\widehat{(x_1,x_2)}$ is the \textit{essential preimage} of
$(x_1,x_2)$ under $(s,t)$, namely the comma category $((s,t)
\downarrow (x_1,x_2))$.
\end{lemma}

\begin{proof}
The argument here is similar to that in Theorem
\ref{thm:groupoidification}, except that we must deal with nontrivial
representations of the $Aut(x_i)$.  That is, when we apply the
Nakayama isomorphism, and the evaluation maps, we cannot use
triviality.  The ``group-averaging'' acting on intertwiners in the
expression we have given is exactly the exterior trace used in the
Nakayama isomorphism.  Here its function is to project a linear map (a
``pulled-back'' intertwiner) onto a space of intertwiners as we push
it along the functor $t$.

In particular, the effect of $\eta_R,s$ on coordinates (i.e. choosing
particular representations $W_i$) is to take an intertwiner $f \in
\opname{hom}[s_1^{\ast}(W_1),t_1^{\ast}(W_2)]$ and produce an
intertwiner at the representations pulled back to $Y$.  The counit
$\eta_L,t$ ``pushes'' this down to an intertwiner in
$\opname{hom}[s_2^{\ast}(W_1),t_2^{\ast}(W_2)]$.  The group averaging
ensures this will be an intertwiner itself.
\end{proof}

\begin{remark}
Using the formula for composition of 2-linear maps and natural
transformations in a general 2-vector space, we can readily see how
horizontal and vertical composition work.

Vertical composition is given by composition of linear maps
component-wise, so we have:
\begin{equation}
\FV(Y' \circ Y)_{([a_1],W_1),([a_2],W_2)}
\end{equation}
with components given by by:
\begin{align}
& \bigoplus_{([x_1],[x_3])} |\widehat{(x_1,x_3)}| \sum_{g \in \opname{Aut}(x_3)} g^{-1} f g \\
\nonumber = & \bigoplus_{([x_1],[x_3])} \sum_{[x_2]} |\widehat{(x_2,x_3)}| \cdot |\widehat{(x_1,x_2)}| \cdot \#\opname{Aut}(x_2) \Bigl{(} \sum_{g \in \opname{Aut}(x_3)} g^{-1} f g \Bigr{)}
\end{align}

This uses the fact that the two group averages each give projections
into spaces of intertwiners, which is redundant, so we omit one,
taking only the order of the group.  We also use that
$\widehat{(x_1,x_3)}$ is a subgroupoid of $Y' \circ Y$.  In fact, it
is a union over all equivalence classes $[x_2]$ in $X_2$ of the
objects in the weak pullback $Y' \circ Y$ based over $[x_2]$, which gives the sum over $[x_2]$ (which performs the matrix multiplication in each component).

The horizontal composite, $\FV(Y' \circ Y) \cong \FV(Y') \circ
\FV(Y)$, on the other hand, involves ``matrix multiplication'' at the
level of composition of 2-linear maps.  The
$(([a_1],W_1),([a_3],W_3))$ component of the product is a linear map
given as a block matrix with one block for each basis 2-vector.  The
blocks consist of the tensor products of the matrices from the
components of $\FV(Y)$ and $\FV(Y')$.  In particular, the $\beta$ isomorphisms from the horizontal composition of source and target induce an isomorphism which acts on intertwiners $\iota \otimes \iota '$ by:
\begin{align}
& \bigl{(} \bigoplus_{\stackrel{([x_1,f_1,x'_1])}{([x_2,f_2,x'_2])}} |\widehat{((x_1,x'_1),(x_2,x'_2))}| \sum_{(g,g') \in \opname{Aut}([x_1],[x'_1])} (g,g')^{-1} \iota \otimes \iota ' (g,g') \bigr{)} \\
\nonumber \cong & \bigoplus_{([a_2],W_2)} \Bigl{(} \bigoplus_{([x_1],[x_2])} |\widehat{(x_1,x_2)}| \sum_{g \in \opname{Aut}(x_1)} g^{-1} \iota g\Bigr{)} \otimes \bigl{(}\bigoplus_{([x'_1],[x'_2])} |\widehat{(x'_1,x'_2)}| \sum_{g' \in \opname{Aut}(x'_1)} (g')^{-1} \iota ' g' \Bigr{)} \\
\nonumber = & \bigoplus_{([a_2],W_2)} \Bigl{(} \bigoplus_{\stackrel{([x_1],[x_2])}{([x'_1],[x'_2])}} |\widehat{(x_1,x_2)}| \cdot |\widehat{(x'_1,x'_2)}| \sum_{(g,g') \in \opname{Aut}(x_1) \times \opname{Aut}(x'_1)} (g,g')^{-1} \iota \otimes \iota ' (g,g')\Bigr{)}
\end{align}

Here, we note that since $Y' \circ Y$ is a weak pullback over $A_2$,
its objects consist of triples $(y,h,y')$, we implicitly have a sum
over $[y,h,y']$ in the groupoid cardinality, which is
$|\widehat{(x_1,x_2)}| \cdot |\widehat{(x'_1,x'_2)}| \cdot
|\opname{Aut}(a_2)$.
\end{remark}

\section{Main Theorem}\label{sec:maintheorem}

Having now described the effect of the functor $\FV$ at each level -
groupoids, spans, and spans of spans---it remains to check that these
really define a 2-functor of the right kind.  We begin by explicitly
laying out what this 2-functor is, then verify the remaining
properties.

\subsection{The 2-Linearization Functor}

We have been defining the maps involved in $\FV$ throughout the last
few sections, so here we collect the full definition in one place.

\begin{definition}\label{def:Vrep} The 2-linearization process $\FV :
\Span(\Gpd) \ra \iiV$ is defined as
follows: \begin{itemize}
\item For an essentially finite groupoid $A$ it assigns:
\begin{equation}\label{eq:VonGpd}
  \FV(A) = \HV{A}
\end{equation}
\item For a span of groupoids:
\begin{equation}
A \lalim^{s} X \ralim^{t} B
\end{equation}
it assigns:
\begin{equation}\label{eq:VonSpan}
  \FV(S) = t_{\ast} \circ s^{\ast}
\end{equation}
\item For a (strictly commuting) span of maps between two spans with the
      same source and target:
\begin{equation}
\xymatrix{
 & X_1 \ar[dl]_{s_1} \ar[dr]^{t_1} & \\
A & Y \ar[u]_{s} \ar[d]^{t} & B  \\
 & X_2  \ar[ul]_{s_2} \ar[ur]^{t_2} & \\
}
\end{equation}
$\FV$ assigns a natural transformation:
\begin{equation}
\FV(Y) = \epsilon_{L,t} \circ \eta_{R,s} : (t_1)_{\ast}s_1^{\ast} \Longrightarrow (t_2)_{\ast}s_2^{\ast}
\end{equation}
(and analogously for weakly commuting spans of maps as in Definition
\ref{def:Von2Span})
\end{itemize}
$\FV$ also associates the following:
\begin{itemize}
\item For each composable pair $X : A_1 \ra A_2$ and $X' : A_2 \ra A_3$, a natural isomorphism
\begin{equation}
 \beta : \FV( X' \circ X ) \ra \FV(X') \circ \FV(X)
\end{equation}, as described in Theorem \ref{thm:FVcomposfunc}.
\item For each object $X \in \Gpd$, the natural transformation
\begin{equation} U_B : 1_{\FV(B)} \stackrel{\sim}{\ra} \FV(1_B) \end{equation}
is the natural transformation induced by the equivalence between
$ {B}$ and $ {1_B}$. 
\end{itemize} 
\end{definition}

Then we have the following:

\begin{theorem}\label{thm:mainthm} The construction given in Definition \ref{def:Vrep} defines a weak 2-functor $\FV : \Span(\Gpd) \ra \iiV$ .
\end{theorem}
\begin{proof}
First, we note that by the result of Lemma \ref{lemma:fgkv}, we know
that $\FV$ assigns a 2-vector space to each object of
$\Span(\Gpd)$.

If $S : B \ra B'$ span of essentially finite groupoids---i.e. a
morphism in $\Span(\Gpd)$, the map $\FV(S)$ defined in Definition
\ref{def:ZGonS} is a linear functor by the result of Theorem
\ref{thm:2mapadjoints}, since it is a composite of two linear maps.
This respects composition of morphisms, as shown in Theorem \ref{thm:FVcomposfunc}, and of 2-morphisms in both horizontal and vertical directions, as shown in Theorems \ref{thm:vertcompos} and \ref{thm:horizcompos}.

Next we need to check that our $\FV$ satisfies the remaining properties of a
weak 2-functor: that the isomorphisms from the weak preservation of
composition and units satisfy the requisite coherence conditions; and
that $\FV$ strictly preserves horizontal and vertical composition of
natural transformations.

The coherence conditions for the compositor morphisms
\begin{equation}
\beta_{S,T} :   \FV(T \circ S) \ra \FV(T) \circ \FV(S)
\end{equation}
and the associator say that these must make the following diagram
commute for all composable triples $(X,X',X'')$:
\begin{equation}\label{eq:wk2fncassoc}
\xy
 (0,20)*+{\FV(X'') \circ \FV(X') \circ \FV(X)}="top";
 (35,4)*+{\FV(X'' \circ X') \circ \FV(X)}="rttop";
 (25,-20)*+{\FV( (X'' \circ X') \circ X )}="rtbot";
 (-25,-20)*+{\FV(X'' \circ (X' \circ X) )}="leftbot";
 (-35,4)*+{\FV(X'') \circ \FV(X' \circ X)}="leftop";
     {\ar^{1 \otimes \beta_{2,1}} "leftop";"top"}
     {\ar_{\beta_{3,2} \otimes 1} "rttop";"top"}
     {\ar^{\beta_{3,21}} "leftbot";"leftop"}
     {\ar_{\beta_{32,1}} "rtbot";"rttop"}
     {\ar_{\FV(\alpha_{X'',X',X})} "leftbot";"rtbot"}
\endxy
\end{equation}

We implicitly assume here a trivial associator for the 2-linear maps
in the expression $\FV(X'') \circ \FV(X') \circ \FV(X)$.  This is
because each 2-linear map is just a composite of functors, so this
composition is associative.  But note that we can similarly assume,
without loss of generality, that the associator $\alpha$ for
composition of spans is trivial.  The composite $X' \circ X$ is a weak
pullback.  This is only defined up to isomorphism, but one candidate
is the comma category for any $x \in A_2$.  Any other candidate is
isomorphic to this one.  But then, the associator
\begin{equation}
\alpha_{X'',X',X} : \FV( X'' \circ (X' \circ X) ) \ra \FV( (X'' \circ X') \circ X )
\end{equation}
is just given by the obvious canonical map between the comma
categories.  In particular, both composites give comma categories
whose objects are determined by choices $(x,f,x',g,x'')$ where $f :
t_1(x) \ra s_2(x')$ and $g : t_2(x') \ra s_3(x'')$, and whose
morphisms are triples of morphisms in $X \times X' \times X''$ making
the appropriate diagrams commute.  However, these comma categories are
defined in terms of pairs, with different parenthesizations.  So
$\alpha_{X'',X',X}$ is the evident isomorphism between these
composites.

So it suffices to show that, up to this identification:
\begin{equation}
(1 \otimes \beta_{X',X}) \circ \beta_{X'',X'\circ X} = (\beta_{X'',X'} \otimes 1) \circ \beta_{X'' \circ X',X}
\end{equation}
The $\beta$ isomorphisms are given by the $\alpha$ up to which the weak pullbacks commute, and so are given by choices of the functions $f \in A_2$ and $g \in A_3$ in the comma categories.  The associator isomorphism induces an corresponding isomorphism between these composite $\beta$ maps by the correspondence between the choices of $f$ and $g$ in the pullback squares on each side of this equation.  So indeed, this is true.

In general, the coherence conditions for the ``unit'' isomorphism
\begin{equation}
U_A : 1_{\FV(A)} \stackrel{\sim}{\ra} \FV(1_A)
\end{equation}
which accomplishes weak preservation of identities, say that it must
make the following commute for any span $X: A_1 \ra A_2$:
\begin{equation}\label{eq:wk2fncunit}
\xymatrix{
\FV(X)  & & \\
 & &  \\
\FV(X) \circ \FV(1_{A_1}) \ar[uu]_{1 \otimes U_{A_1}} &  & \FV(X \circ 1_{A_1})\ar[uull]_{\FV(r_X)} \ar[ll]_{\beta_{X,1_{A_1}}}
}
\end{equation}
where $r_{A_1}$ is the right unitor for $A_1$.  There is also the symmetric
condition for the left unitor.

We notice that, as with $\FV(1_{A_1})$, $\FV(r_{A_1})$ is equivalent to the
identity.  The map $r_X : X \circ 1_{A_1} \ra X$ is the canonical
isomorphism taking composition of $X$ with an identity span to $X$
which is just a projection from a comma category.  Since $X \circ 1_{A_1}$
and $X$ are thus isomorphic, .

So the condition amounts to the fact that $\beta_{X,1_{A_1}} : \FV(X
\circ 1_{A_1}) \ra \FV(X) \circ \FV(1_{A_1}) = \FV(X) $ is equivalent
to the identity in such a way that (\ref{eq:wk2fncunit}) commutes.
But this is immediate since this $\beta$ map is being applied to an
identity span.
\end{proof}

\section{Acknowledgements}

The author would like to recognize the invaluable assistance of John
Baez, whose advice and guidance made possible the Ph.D. thesis work
which led to this project; to acknowledge the useful discussions with
and help of James Dolan, Derek Wise, Alex Hoffnung, Jamie Vicary, and
Dan Christensen; and to the editors and referees whose copious
suggestions greatly improved the material and its presentation.

\appendix

\section{Weak Preservation of Composition}\label{ax:composproof}

In this appendix, we give some background to the definition of composition of spans of groupoids, namely comma categories.  We also give a note on a key element of the proof of Theorem
\ref{thm:FVcomposfunc}, which states that the putative 2-functor $\FV$
weakly preserves this composition.  We rely on the fact that a pullback
square of groupoids gives rise to a square of 2-linear maps, which
satisfies the Beck-Chevalley condition.  We discuss this here as well.

\subsection{Background on Comma Categories}\label{sec:commacat}

We now recall some facts about comma categories, which play a role in
our construction of our 2-functor $\FV$ in the composition of spans of
groupoids, via weak pullback.

\begin{definition} Given a diagram of categories $\catname{A}
\ralim^{F} \catname{C} \lalim^{G} \catname{B}$.  Then an object in the
\textbf{comma category} $(F \downarrow G)$ consists of a triple
$(a,f,b)$, where $a \in \catname{A}$ and $b \in \catname{B}$ are
objects, and $f : F(a) \ra G(b)$ is a morphism in $\catname{C}$.  A
morphism in $(F \downarrow G)$ from $(a_1,f_1,b_1)$ to $(a_2,f_2,b_2)$
consists of a pair of morphisms $(h,k) \in \catname{A}
\times \catname{B}$ making the square
\begin{equation}\label{eq:commacatmor}
\xymatrix{
F(a_1) \ar[r]^{f_1} \ar[d]_{F(h)} & G(b_2) \ar[d]^{G(k)} \\ F(a_2)
\ar[r]_{f_2} & G(b_2) }
\end{equation}
commute.
\end{definition}

\begin{remark}Note that in a weak pullback, the morphisms $f$ would be
required to be an \textit{isomorphism}, but when we are talking about
a weak pullback of groupoids, these conditions are the same.\end{remark}

The comma category has projection functors which complete the (weak)
pullback square for the two projections:
\begin{equation}
\xymatrix{   & (F \downarrow G) \ar_{P_A}[ld] \ar^{P_B}[rd] & \\
   \catname{A}\ar_{F}[rd] \ar@{=>}[rr]^{\alpha}_{\sim} & & \catname{B}\ar^{G}[ld] \\
   & \catname{C} & \\
}
\end{equation}
such that $(F \downarrow G)$ is a universal object (in $\Cat$) with
maps into $\catname{A}$ and $\catname{B}$ making the resulting square
commute up to a natural isomorphism $\alpha$.  This satisfies the
universal condition that, given any other category $\catname{D}$ with
maps to $\catname{A}$ and $\catname{B}$, there is an equivalence
between $[\catname{D},\catname{C}]$ and the comma category
$(P_A^{\ast},P_B^{\ast})$ (where $P_S*$ and $P_T*$ are the functors
from $\catname{D}$ to $\catname{B}$ which factor through $P_S$ and
$P_T$ respectively).  This is the weak form of the universal property
of a pullback.

So suppose we restrict to the case of a weak pullback of groupoids.
This is equivalent to the situation where $\catname{A}$, $\catname{B}$
and $\catname{C}$ are skeletal - that is, each is just a disjoint
union of groups.  Then the set of objects of $(F \downarrow G)$ is a
disjoint union over all the morphisms of $\catname{C}$ (which are all
of the form $g: x \ra x$ for some object $x$) of all the pairs of
objects $a \in \catname{A}$ and $b \in \catname{B}$ with $g: F(a) \ra
G(b)$.  In particular, since we assume $\catname{C}$ is skeletal, this
means $F(a) = G(b)$, though there will be an instance of this pair in
$(F \downarrow G)$ for each $g$ in the group of morphisms on this
object $F(a) = G(b)$.

So as the set of objects in $(F \downarrow G)$ we have a disjoint
union of products of sets---for each $c \in \catname{C}$, we get
$|Aut(c)|$ copies of $F^{-1}(c) \times G^{-1}(c)$.  The set of
morphisms is just the collection of commuting squares as in
(\ref{eq:commacatmor}) above.

Note that if we choose a particular $c$ and $g: c \ra c$, and choose
objects $a$, $b$ with $F(a)=c$, $G(b)=c$, and if $H=\opname{Aut}(a)$,
$K=\opname{Aut}(b)$ and $M=\opname{Aut}(c)$, then the group of
automorphisms of $(a,g,b) \in (F \downarrow G)$ is isomorphic to the
fibred product $H \times_M K$.  In particular, it is a subgroup of the
product group $H \times K$ consisting of only those pairs $(h,k)$ with
$F(h) g = g G(k)$, or just $F(h) = g G(k) g^{-1}$.  We can call it $H
\times_M K$, keeping in mind that this fibred product depends on $g$.
Clearly, the group of automorphisms of two isomorphic objects in $(F
\downarrow G)$ are isomorphic groups.

Now, as we saw when discussing comma squares, the objects of the weak
pullback ${X' \circ X}$ consist of pairs of objects, $x \in
X$, and $x' \in X'$, together with a morphism in $A_2$, $g
: t(x) \ra s'(x')$.  The morphisms from $(x_1,g_1,{x'}_1)$ to
$(x_2,g_2,{x'}_2)$ in the weak pullback are pairs of morphisms, $(h,k) \in
{X} \times  {X'}$, making the square
\begin{equation}
\xymatrix{
t(x_1) \ar[r]^{g_1} \ar[d]_{t(h)} & s'({x'}_2) \ar[d]^{s'(k)} \\
t(x_2) \ar[r]_{g_2} & s'({x'}_2)
}
\end{equation}
commute.

We may assume that the groupoids we begin with are skeletal---if not,
we replace the groupoid with its skeleton, so the objects are just
isomorphism classes of the original objects.  Then recall from Section
\ref{sec:span2lin} that in this weak pullback the set of objects in
${X' \circ X}$ is a disjoint union of products of sets - for each $a
\in {A_2}$, we get $|Aut(a)|$ copies of ${t}^{-1}(a) \times
(s')^{-1}(a)$.

\subsection{The Beck Condition}

\begin{remark} The isomorphism $\alpha$ in the weak pullback square
(\ref{xy:wkpb}) gave rise to a natural isomorphism:
\begin{equation}
\alpha^{\ast} :  T^{\ast} \circ  (s')^{\ast}   \ra  S^{\ast} \circ t^{\ast}
\end{equation}
Given an object in the composite $X' \circ X$, $\alpha$
gives an isomorphism of the two restrictions to $A_2$, through $X$ and $X'$.

What we proved is that the other square---the ``mate'' under the
adjunctions, also has a natural isomorphism (``vertically'' across the
square), namely that there exists:
\begin{equation}
\beta_{X,X'} : T_{\ast} \circ S^{\ast} \ra (s')^{\ast} \circ t_{\ast}
\end{equation}

In fact, these are related by the units for both pairs of adjoint functors:
\begin{equation}
\eta_{R,T} : 1_{\FV(X'\circ X)} \ra  T_{\ast} \circ T^{\ast}
\end{equation}
and
\begin{equation}
\eta_{R,t} : 1_{\FV(X)} \ra  t_{\ast} \circ t^{\ast}
\end{equation}

So the desired ``vertical'' natural transformation across the square
(\ref{eq:adjunctionmatesquare}) is determined  by the condition that it complete the following
square of natural transformations to make it commute:
\begin{equation}\label{eq:nearlybeck2}
\xymatrix{
T_{\ast} \circ S^{\ast} \ar@{=>}[rr]^(.4){1 \otimes \eta_{R,t}} \ar@{==>}[d]^{\beta_{X,X'}} & & T_{\ast} \circ S^{\ast} \circ t^{\ast} \circ t_{\ast} \ar@{=>}[d]^{1 \otimes (\alpha^{\ast})^{-1} \otimes 1} \\ 
(s')^{\ast} \circ t_{\ast} \ar@{=>}[rr]^(.4){1 \otimes \eta_{R,T}} & & T_{\ast} \circ T^{\ast} \circ (s')^{\ast} \circ t_{\ast} \\ 
}
\end{equation}

The crucial element of this is the fact that the (weak) pullback
square for the groupoids in the middle of the composition diagram
gives rise to a (weak) pullback square of $\V$-presheaf categories.
This is shown by Ross Street \cite{street2}.  This is the
\textit{Beck-Chevalley} (BC) condition, which is discussed by
B\'enabou and Streicher \cite{BCES}, MacLane and Moerdijk
\cite{macmoer}, and by Dawson, Par\'e and Pronk
\cite{unispan}.
\end{remark}

\bibliography{grpd-2vect}

\begin{thebibliography}{10}

\bibitem{hdaII}
{\sc Baez, J.}
\newblock Higher-dimensional algebra {II}: 2-hilbert spaces.
\newblock {\em Adv. Math. 127\/} (1997), 125--189.

\bibitem{hdaVII}
{\sc Baez, J.}
\newblock Higher-dimensional algebra {VII}: {G}roupoidification.
\newblock {\tt{http://math.ucr.edu/home/baez/hda7.pdf}}.

\bibitem{BC}
{\sc Baez, J., and Crans, A.}
\newblock Higher-dimensional algebra {VI}: {L}ie 2-algebras.
\newblock {\em Theory and Applications of Categories 12\/} (2004), 492--528.

\bibitem{finfeyn}
{\sc Baez, J., and Dolan, J.}
\newblock From finite sets to {F}eynman diagrams.
\newblock In {\em Mathematics Unlimited - 2001 And Beyond\/} (2001),
  B.~Engquist and W.~Schmid, Eds., Springer Verlag.
\newblock Preprint at {\tt{arXiv:math.QA/0004133}}.

\bibitem{BCES}
{\sc B\'enabou, J., and Streicher, T.}
\newblock {B}eck-{C}hevalley condition and exact squares.
\newblock Unpublished.

\bibitem{benson}
{\sc Benson, D.}
\newblock {\em Representations and Cohomolgy I: Basic representation theory of
  finite groups and associative algebras}.
\newblock Cambridge studies in advanced mathematics. Cambridge University
  Press.

\bibitem{unispan}
{\sc Dawson, R. J.~M., Par\'e, R., and Pronk, D.~A.}
\newblock Universal properties of span.
\newblock {\em Theory and Applications of Categories 13}, 4 (2004), 61--85.

\bibitem{elgueta}
{\sc Elgueta, J.}
\newblock Generalized 2-vector spaces and general linear 2-groups.
\newblock {\em J. Pure Appl. Alg. 212\/} (2008), 2067--2091.

\bibitem{freyd}
{\sc Freyd, P.}
\newblock {\em {A}belian Categories: An Introduction to the Theory of
  Functors}.
\newblock Harper \& Row, 1964.

\bibitem{KV}
{\sc Kapranov, M., and Voevodsky, V.}
\newblock 2-categories and {Z}amolodchikov tetrahedron equations.
\newblock {\em Proc. Symp. Pure Math 56 Part 2\/} (1994), 177--260.

\bibitem{laudaambidjunction}
{\sc Lauda, A.~D.}
\newblock Frobenius algebras and ambidextrous adjunctions.
\newblock {\em Theory and Applications of Categories 16}, 4 (2006), 84--122.

\bibitem{lawverethesis}
{\sc Lawvere, W.~F.}
\newblock {\em Functorial Semantics of Algebraic Theories and Some Algebraic
  Problems in the Context of Functorial Semantics of Algebraic Theories}.
\newblock PhD thesis, Columbia University, 1963.
\newblock Reprinted in Theory and Applications of Categories, 2004.

\bibitem{maclane}
{\sc MacLane, S.}
\newblock {\em Categories for the Working Mathematician}.
\newblock No.~5 in Graduate Texts in Mathematics. Springer, 1971.

\bibitem{macmoer}
{\sc MacLane, S., and Moerdijk, I.}
\newblock {\em Sheaves in Geometry and Logic: A First Introduction to Topos
  Theory}.
\newblock Universitext. Springer Verlag, 1992.

\bibitem{mortonthesis}
{\sc Morton, J.~C.}
\newblock {\em Extended TQFT's and Quantum Gravity}.
\newblock PhD thesis, University of California, Riverside, 2007.
\newblock {\tt{arxiv:math/0710.0032}}.

\bibitem{elango}
{\sc Panchadcharam, E.}
\newblock {\em Categories of Mackey Functors}.
\newblock PhD thesis, Macquarie University, 2006.

\bibitem{sternberg}
{\sc Sternberg, S.}
\newblock {\em Group theory and physics}.
\newblock Cambridge University Press, 1994.

\bibitem{street2}
{\sc Street, R.}
\newblock Enriched categories and cohomology.
\newblock {\em Quaestiones Mathematicae 6\/} (1983), 265--283.
\newblock Reprints in Theory and Applications of Categories.

\end{thebibliography}
\bibliographystyle{acm}

\end{document}